\documentclass[11pt]{article}
\usepackage{enumitem}
\usepackage{graphicx,epstopdf,bm}
\usepackage{geometry}
\makeatletter
\newcommand{\singlespacing}{\let\CS=\@currsize\renewcommand{\baselinestretch}{1}\tiny\CS}
\newcommand{\oneandahalfspacing}{\let\CS=\@currsize\renewcommand{\baselinestretch}{1.25}\tiny\CS}
\newcommand{\doublespacing}{\let\CS=\@currsize\renewcommand{\baselinestretch}{1.35}\tiny\CS}

\newtheorem{theorem}{Theorem}[section]

\newtheorem{rule-def}[theorem]{Rule}

\begin{document}
\newcommand{\la}{\lambda}
\newcommand{\si}{\sigma}
\newcommand{\ol}{1-\lambda}
\newcommand{\bphi}{\mbox{\boldmath$\phi$}}
\newcommand{\be}{\begin{equation}}
\newcommand{\ee}{\end{equation}}
\newcommand{\bea}{\begin{eqnarray}}
\newcommand{\eea}{\end{eqnarray}}
\newcommand{\nn}{\nonumber}
\newcommand{\mc}{\multicolumn}
\newcommand{\bee}{\begin{eqnarray*}}
\newcommand{\eee}{\end{eqnarray*}}
\newcommand{\lb}{\label}
\newcommand{\D}{${\bf{D~}}$}
\newcommand{\C}{${\bf{C~}}$}
\newcommand{\rr}{\rule{2mm}{2mm}}
\newcommand{\spc}{\hspace{0.5cm}}
\newcommand{\nii}{\noindent}
\newcommand{\ii}{\indent}
\newcommand{\0}{${\bf 0}$}
\newcommand{\bpi}{\mbox{\boldmath$\pi$}}
\newcommand{\w}{\omega}
\newcommand{\bw}{\mbox{\boldmath$\omega$}}
\newcommand{\bsi}{\mbox{\boldmath$\psi$}}
\newcommand{\bet}{\mbox{\boldmath$\beta$}}
\newcommand{\al}{\mbox{\boldmath$\alpha$}}
\newcommand{\n}{\mbox{\boldmath$\nu$}}
\newcommand{\bPi}{\mbox{\boldmath$\Pi$}}
\newcommand{\bpsi}{\mbox{\boldmath$\psi$}}
\newcommand{\sbpi}{\overline{\mbox{\boldmath$\pi$}}}
\newcommand{\Phib}{\mbox{\boldmath$\Phi$}}
\newcommand{\alphab}{\mbox{\boldmath$\alpha$}}
\doublespacing
\title{ Steady-state analysis of single exponential vacation
in a  $PH/MSP/1/\infty$ queue using roots}
\author{ A. D. Banik$^{\mbox{a,}}$ \thanks{Corresponding
author's Tel. +91-6742576071. Fax. +91-6742301983. \newline E-mail
addresses: banikad@gmail.com;
adattabanik@iitbbs.ac.in (A. D. Banik)~,~chaudhry-ml@rmc.ca (M. L. Chaudhry)~,~florin.avram@univ-Pau.fr (F. Avram)}~,~ M. L. Chaudhry $^{\mbox{b}}$~,~Florin Avram$^{\mbox{c}}$\\
$^{\mbox{a}}${\it School of Basic Sciences, Indian Institute of Technology Bhubaneswar,}\\
{\it Permanent Campus Argul, Jatni, Khurda-752 050, Odisha, India}\\$^{\mbox{b}}$ {\it Department of Mathematics and Computer Science}\\{\it Royal Military College of Canada}\\{\it P.O. Box 17000, STN Forces, Kingston Ont., Canada K7K 7B4}\\ $^{\mbox{c}}$ {\it Laboratoire de Math\'{e}matiques Appliqu\'{e}es}\\{\it Universit\'{e} de Pau, France.}
}
\maketitle
\begin{abstract}
We consider an infinite-buffer single-server queue where   inter-arrival
times are phase-type ($PH$),  the service
is provided according to Markovian service
process $(MSP)$, and the server may take  single, exponentially
distributed  vacations when the queue is empty.  The proposed analysis is based on roots of
the associated characteristic equation of the vector-generating
function (VGF) of system-length distribution at a pre-arrival epoch.
Also, we obtain the steady-state system-length distribution at an
arbitrary epoch along with
some important performance measures such as the mean number of customers in the system and the mean system sojourn time
of a customer. Later, we have established heavy- and
light-traffic approximations as well as an approximation for the
tail probabilities at pre-arrival epoch based on one root of the
characteristic equation. At the end, we present numerical results in
the form of tables to
show the effect of model parameters on the performance measures.    \vspace{0.3cm}
\end{abstract}
\tableofcontents
\section{\Large{\bf Introduction}}\lb{id1}
In recent times, queueing models with non-renewal arrival and
service processes have been used to model networks of complex
computer and communication systems. Traditional queueing analysis
using Poisson processes is not powerful enough to capture the
correlated nature of arrival (service) processes. The performance
analysis of correlated type of arrival processes may be done through
some analytically tractable arrival process viz., Markovian arrival
process ( $MAP),$ see Lucantoni et al. ${\hbox{\cite{nu4}}}.$ The
$MAP$ has the property of both time varying arrival rates and
correlation between inter-arrival times. To consider batch arrivals
of variable capacity, Lucantoni \cite{abl2} introduced batch
Markovian arrival process ( $BMAP)$. The processes $MAP$ and $BMAP$ are
convenient representations of a versatile Markovian point
process, see Neuts \cite{x_2} and Ramaswami \cite{r80}. Like the
$MAP$, Markovian service process $(MSP)$ is a versatile
service process which can capture the correlation among the
successive service times. Several other service processes, e.g.,
Poisson process, Markov modulated Poisson process ($MMPP$) and
phase-type ($PH$) renewal process can be considered as special cases
of $MSP$. For details of $MSP$, the readers are referred to
Bocharov \cite{nebb11} and Albores and Tajonar \cite{nebb2}. The
analysis of finite-buffer $G/MSP/1/r~(r\leq \infty)$ queue has
been performed by Bocharov et al. \cite{nebb1}. The same queueing
system with multiple servers such as $GI/MSP/c/r$ has been
analyzed by Albores and Tajonar \cite{nebb2}. Gupta and Banik
\cite{bub122} analyzed $GI/MSP/1$ queue with finite- as well as
infinite-buffer capacity using a combination of embedded Markov chain and supplementary variable method. \par  During the last two
decades, queueing systems with vacations have been studied
extensively. For more details on this topic, the readers are referred to
a recent book by Tian and Zhang \cite{tz06} and references
therein. An extensive amount of literature is available on infinite-
and finite-buffer $M/G/1$- and $GI/M/1$-type queueing models with
multiple vacations, see first few chapters of \cite{tz06}, Karaesmen
and Gupta \cite{abl33} and Tian et al. \cite{bu3333}. However,
limited studies have been done on $GI/M/1$ queue with single
vacation, see Chapter 4 of \cite{tz06}. In the past few years, there
is a growing trend to analyze queueing models with renewal or
non-renewal arrival and service processes with server vacation, see,
e.g., Lucantoni et al. ${\hbox{\cite{nu4}}}$ and Shin and Pearce
\cite{rev_3}. The analysis of phase-type server vacation for the
case of $GI/M/1$ queue has been carried out by Chen et al.
\cite{chenetal09}. Baba \cite{baba10} analyzes $M/PH/1$ queue
where the server is allowed to take working vacations as well as
vacation interruptions. Samanta \cite{s09} discussed a discrete-time
$GI/Geo/1$ queue with single geometric vacation time.
Recenly, Chaudhry et al. \cite{csa12,cbp15} discussed $GI/MSP/1/\infty$ queues with single and batch arrivals using the roots method, respectively.

\par In this paper, we carry out the analytic analysis of the
$PH/MSP/1/\infty$ queue with exponential single vacation through
the calculation of roots of the denominator of the underlying vector
generating function of the steady-state probabilities at pre-arrival
epoch. In this connection, the readers are referred to Chaudhry et
al. \cite{csa12,cgg10,cbp15}, Tijms \cite{t03} and Chaudhry et al.
\cite{chm90} who have used the roots method. The roots can be easily
found using one of the several commercially available packages such
as Maple and Mathematica. The algorithm for finding such roots is
available in some papers, e.g., see Chaudhry et al. \cite{chm90}.
The purpose of studying this queueing model using roots is that we
obtain computationally simple and analytically closed form solution
to the infinite-buffer $PH/MSP/1$ queue with the vacation time
following exponential distribution. It may be remarked here that the
matrix-geometric method (MGM) uses iterative procedure to get
steady-state probabilities at the pre-arrival epochs. Further, it is
well known that for the case of the MGM it is required to
solve the non-linear matrix equation with the dimension of each
matrix in this equation being the number of service-phases involved in
a $PH/MSP/1$ queue. In the case of the roots method, we do not have
to investigate the structure of the transition probability matrices
(TPM) at the embedded pre-arrival epochs. It may be mentioned here
that the basic idea of correlated service was first introduced by
Chaudhry \cite{c72}. Further, it may be remarked here that the analysis of the
infinite-buffer queues with renewal input and
exponential service time under exponential server vacation(s) has been carried out by Tian and Zhang
\cite{tz06}, see Chapter 4. The queueing model that we are going to
consider has non-renewal service ($MSP$) and exponential single
vacation time. In addition, we discuss several other quantitative
measures such as system-length distribution at a post-departure
epoch and expected busy and idle periods. Later, we have established heavy- and
light-traffic approximations as well as an approximation for the
tail probabilities at pre-arrival epoch based on one root of the
characteristic equation. Finally, some numerical
results have been presented which may help researchers/practitioners
to tally their results with those of ours.


\section{\bf Description of the model}\lb{id2}
Let us consider a single-server infinite-buffer queueing system with
the server's single vacation. The inter-arrival time of customers,
the service time of a customer and the vacation time of the server
are represented by the generic random variables (r.v.'s) $A$,~$S$ and $V$,
respectively. Let $F_X(x)$ denote the distribution
function (D. F.) of the random variable $X$ with $f_X(x)$ and
$f_X^{*}(s)$ the corresponding probability density function
(p.d.f.) and Laplace-Stieltjes transform (LST), respectively. The
inter-arrival time $A$ is assumed to have a general distribution with
p.d.f. $f_A(x)$, D. F. $F_A(x)$ and LST $f_A^{*}(s)$.\\
{\bf Arrivals.} The inter-arrival times are assumed to be independent and
identically-distributed (i.i.d.) random variables and they are
independent of the service process as well as vacation time. The inter-arrival time distribution $PH$ is an important special case of general distribution
as  the distribution possesses nice vector and matrix form
representation. Several probability distributions such as Earlang,
hyper-exponential, generalized Earlang, Coxian etc. can be treated
as special cases of $PH$-distribution. It may be noted here that $PH$-distribution is a special case of general distribution. If the
inter-arrival times follow $PH$-type distribution with
irreducible representation $({\bm \alpha}, {\bm  T})$, where ${\bm \alpha}$ \& ${\bm 
T}$ are a vector and a matrix of dimension $1 \times \eta$ and
$\eta\times \eta$, respectively, the p.d.f. and D.F. of
inter-arrival times are given by \bea F_{A}(x) & = & 1 - {\bm \alpha} \
e^{{\bm  T}x}{\bm  e}_{\eta},\quad \mbox{for} \ x \ge 0,\lb{rem_1}\\
\mbox{and} \quad f_{A}(x) & = & - {\bm \alpha} \ e^{{\bm  T}x}
{\bm  T}{\bm  e}_{\eta}={\bm \alpha} \ e^{{\bm  T}x}
{{\bm  T}}^{0}, \quad  \mbox{for} \ x > 0, \lb{DD114_x} \eea where ${{\bm  T}}^{0}$ is a
non-negative vector and satisfies ${\bm  T}{\bm 
e}_{\eta}+{{\bm  T}}^{0}={\bm  0}$ and ${\bm  e}_{\eta}$ is an $\eta
\times 1$ vector with all its elements equal to 1. Throughout the
paper we write a subscript as the dimension of the column vector
${\bm  e}$ and sometimes we write ${\bm  e}$ by dropping its
subscript. The mean inter-arrival time during a normal busy period is
given by \bea \frac{1}{\lambda}= {\bm \alpha} \int_0^{\infty}xe^{{\bm  T}x} \
dx (-{\bm  T}){\bm  e}_{\eta} &=& - {\bm \alpha} ({\bm  T})^{-1}{\bm 
e}_{\eta}.~~ \lb{viss_1}\eea

\nii {\bf Services.} The customers are served singly according to the continuous-time Markovian
service process ($MSP$) with matrix representation $({\bm 
L}_0,{\bm  L}_1)$. The $MSP$ is a generalization of the Poisson
process where the services are governed by an underlying $m$-state
Markov chain. For more details on $MSP$, the readers are referred to recent papers by Chaudhry et al. \cite{csa12,cbp15}. Let $N(t)$ denote the number of customers
served in $t$ units of time and $J(t)$ the state of the underlying
Markov chain at time $t$ with its state space $\{i : 1\leq i\leq
m\}$. Then $\{N(t),J(t)\}$ is a two-dimensional Markov process with
state space $\{(\ell,i) : \ell\geq 0, 1\leq i\leq m\}$. Average service rate of customers $\mu^{\star}$ (the so
called fundamental service rate) of the stationary $MSP$ is
given by $\mu^{\star}=\overline{\bpi}{\bm  L}_1{\bm  e}$, where
$\overline{\bpi}=[\overline{\pi}_1,\overline{\pi}_2,
\ldots,\overline{\pi}_m]$ with $\overline{\pi}_j$ denoting the
steady-state probability of servicing a customer in phase
$j~(1\leq j\leq m)$. The stationary probability row-vector
$\overline{\bpi}$ can be calculated from $\overline{\bpi}{\bm 
L}={\bm  0}$ with $\overline{\bpi}{\bm  e}=1$, where ${\bm  L}={\bm  L}_0+{\bm  L}_1$.
 The customers are served singly according
to a $MSP$ with steady-state mean service time $1/\mu^*$.

\par Now, let us define $\{{\bm  P}(n, t): n\geq 0,
t\geq 0\}$ as the $m\times m$ matrix whose $(i,j)$th element is the
conditional probability defined as \bee
P_{i,j}(n,t)=Pr\{N(t)=n,J(t)=j|N(0)=0,J(0)=i\},\quad 1\leq i,j\leq
m. \eee Let ${\bm  P}(n, t), n\geq
0, t\geq 0$ be the $m\times m$ matrices whose elements are $P_{i,j}(n,t)$. Then using Chaudhry et al. \cite{csa12,cbp15}, it may be derived that \bea {\bm 
	P}^{\ast}(z,t)=e^{{\bm  L}(z)t},\quad |z|\leq 1,~t\geq
0,\label{cp3}\eea \noindent where ${\bm  L}(z)={\bm  L}_0+{\bm  L}_1z$ and   $ {\bm 
	P}^{\ast}(z,t)=\sum\limits_{n=0}^{\infty}{\bm  P}(n,t)z^n,\quad
|z|\leq 1.$\\

{\bf  Vacations.} The server is allowed to take a single vacation whenever the system
becomes empty. On return from a vacation if the server finds the
system nonempty he will serve the customers present in the queue,
otherwise the server waits for a customer to arrive and the system
continues in this manner. For an exponential single vacation time
represented by the r.v. $V$, the LST, p.d.f. and D.F. are given as
follows: \bea f_{V}^*(s)=\frac{\gamma}{\gamma+s},~~ ~f_{V}(x)=\gamma
e^{-\gamma x},~~ F_{V}(x)=1-e^{-\gamma x}. \lb{nb_8}\eea where
$1/\gamma~(>0)$ is assumed as the mean vacation time. The Vacation
times are independent of the arrival as well as of the service
processes. The traffic intensity is given by $\rho= \la
E(S)=\la/\mu^{*}$ which is also independent of the vacation process.
\section{The vector generating function of the number of customers served during an inter-arrival and other related probability matrices}

 Let ${\bm  S}_n~(n\geq 0)$ denote the matrix of order $m\times
m$ whose $(i,j)$th element represents the conditional probability
that during an inter-arrival period $n$ customers are served and the
service process passes to phase $j$, provided at the initial instant
of the previous arrival epoch there were at least $n$ customers in
the system and the service process was in phase $i$. Then \bea {\bm 
S}_n=\int_0^{\infty}{\bm  P}(n,t)dF_A(t),~n\geq 0.\label{cs1}\eea If
${\bm  S}(z)$ is the matrix-generating function of ${\bm  S}_n$, where
$S_{i,j}(z)~(1\leq i,j\leq m)$ are the elements of ${\bm  S}(z)$,
then, using (\ref{cs1}) and (\ref{cp3}), we get\bea {\bm 
S}(z)&=&\sum\limits_{n=0}^{\infty}{\bm 
S}_nz^n=\int_0^{\infty}\sum\limits_{n=0}^{\infty}{\bm 
P}(n,t)z^ndF_A(t)\nonumber\\&=&\int_0^{\infty}{\bm 
P}^{\ast}(z,t)dF_A(t)=\int_0^{\infty}e^{{\bm 
L}(z)t}f_A(t)dt=f_A^*(-{\bm 
L}(z)).\label{sz1}\eea The evaluation of the matrices ${\bm 
S}_n$ can be carried out along the lines proposed by Lucantoni
\cite{abl2}. For the sake of completeness, we have given the
procedure of obtaining ${\bm  S}_n$, see Lucantoni \cite{abl2}.
One may note that the computation of ${\bm 
S}(z)$ using Equation (\ref{sz1}) may be cumbersome. However, the following scheme may be efficient and is given by
\bea {\bm 
S}(z)&=& \lim_{N \rightarrow \infty}\sum\limits_{n=0}^{N}{\bm 
S}_nz^n, \eea
where ${\bm 
S}_n$ may be obtained as proposed in Chaudhry et al. \cite{cgg10}. 

We further introduce a few more notations which are required for the
rest of the analysis of the queueing model under consideration. Now
from renewal theory of semi-Markov process, if we let $\widehat{A}$
and $\widetilde{A}$ denote the remaining and elapsed times of an
inter-arrival time, respectively, then \bea
F_{\widehat{A}}(x)=F_{\widetilde{A}}(x)=\int_0^x\la(1-F_A(y)) \ dy,
\lb{gos_3}\eea which we use while deriving the
expression for ${\bm  \Omega}_n$ in Equation (\ref{27_s_c1x}). Similar
to the case of inter-arrival time, if we let $\widehat{V}$ denote
the remaining vacation time, then \bea F_{\widehat{V}}(x)&=&
1-e^{-\gamma x}, \Big[\mbox{Using}~(\ref{gos_3})\Big] \lb{gos_3_x}\eea and \bea
f_{\widehat{V}}(x)&=& \gamma e^{-\gamma x}. \lb{vis_2}\eea

\par As above, we introduce the matrices ${\bm  \Omega}_n~(n\geq 0)$
of order $m\times m$ whose $(i,j)$th element represents the limiting
probability that $n$ customers are served during an elapsed
inter-arrival time of the arrival process with the service process
being in phase $j$, given that there were at least $(n+1)$ customers
in the system with the service process being in phase $i$ at the
beginning of the inter-arrival period. Then, from Markov renewal
theory as given in Chaudhry and Templeton \cite[p. 74-77]{chau}, we
have \bea {\bm  \Omega}_n &=& \lambda \int_0^{\infty}{\bm 
P}(n,x)(1-F_A(x)) \ dx, \quad n \geq 0. \lb{27_s_c1x} \eea \nii

\nii The matrices ${\bm  \Omega}_n$ can be expressed in terms of the
matrices ${\bm  S}_n$ and their relationship discussed in \cite{cbp15} is as follows: \bea {\bm  S}_n &=&\delta_{n,0}{\bm  I}_m +
\frac{1}{\lambda}{\bm  \Omega}_n{\bm  L}_0 +\frac{1}{\lambda}{\bm  \Omega}_{n-1}{\bm  L}_{1}.{\it 1_{\{n\geq 1\}}} , \quad n\geq 0, \lb{cobis_1} \eea
where ${\it 1_{\{n\geq 1\}}}$ is an indicator function and takes value 1 if the the condition $n\geq 1$ is satisfied, otherwise it takes value 0. 
%
\par Further, let $\widetilde{P}_{ij}(n,t)$ be the conditional probability that at least $n$
customers are served in $(0,t]$ and the service process is in phase
$j$ at the end of the $n$th service completion, given that there
were $n$ customers in the system and the service process was in
phase $i$ at time $t=0$. The probabilities
$\widetilde{P}_{ij}(n,t),~n\geq 1,~t\geq 0,$ then satisfy the
equations \bee \widetilde{P}_{ij}(n,t+\Delta
t)&=&\widetilde{P}_{ij}(n,t)+\sum\limits_{k=1}^mP_{ik}(n-1,t)[L_1]_{kj}\Delta
t + o(\Delta t),\eee with the initial condition
$\widetilde{P}_{ij}(n,0)=0,~n\geq 1.$ Rearranging the terms and
taking the limit as $\Delta t\to 0$, it reduces to \bee
\frac{d}{dt}\widetilde{P}_{ij}(n,t)&=&\sum\limits_{k=1}^mP_{ik}(n-1,t)[L_1]_{kj},\quad
n\geq 1\eee for $t\geq 0,~1\leq i,j\leq m$, with the initial
conditions $\widetilde{P}_{ij}(n,0)=0$. This system may be written
in matrix notation as \bea \frac{d}{dt}\widetilde{{\bm 
P}}(n,t)&=&{\bm  P}(n-1,t){\bm  L}_1,\quad n\geq 1,\label{ptilta}\eea
with $\widetilde{{\bm  P}}(n,0)={\bm  0},~~n\geq 1$.\\
\nii Let $\omega$ denote the probability that $\widehat{V}$ exceeds
an inter-arrival time $A$, then \bea \omega &=& \int_0^{\infty}
Pr(A<\widehat{V}|A=x).f_{A}(x)
\ dx \nn \\ &=& \int_0^{\infty} Pr(\widehat{V}>x).f_{A}(x) \ dx \nn \\
&=& \int_0^{\infty}(1-F_{\widehat{V}}(x))f_{A}(x) \ dx =f_A^{*}(\gamma). \lb{vis_3_x_1}\eea
 Similarly, if we let $\tau$ denote the probability that
$\widehat{V}$ exceeds $\widehat{A}$, then \bea \tau &=&
\int_0^{\infty} (1-F_{\widehat{V}}(x))f_{\widehat{A}}(x) \ dx =f_{\widehat{A}}^{*}(\gamma).
\lb{vis_4_x}\eea
\nii {\bf  Remark 2.1:}~$\omega$ and $\tau$ may be derived in slightly different way. In the following we present a slightly different derivation for $\omega$ and $\tau$ may be done similarly.  \bea \omega &=& \int_0^{\infty}
Pr(A<\widehat{V}|\widehat{V}=x).f_{\widehat{V}}(x)
\ dx \nn \\ &=& \int_0^{\infty} Pr(A<x).f_{\widehat{V}}(x) \ dx \nn \\
&=& \int_0^{\infty} F_A(x)f_{\widehat{V}}(x) \ dx. \lb{vis_3_x}\eea

\par One of the frequently used inter-arrival time is
phase-type renewal process which also serves as a special case of
several other inter-arrival time distributions and is well-known in
the literature. Therefore, we state the above formulae
(\ref{vis_3_x}) and (\ref{vis_4_x}) for the case of phase-type
inter-arrival time by the following theorem.
\begin{theorem} If inter-arrival time follows a $PH$-distribution with
irreducible representation $({\bm \alpha},{\bm  T})$, where ${\bm \alpha}$ and ${\bm 
T}$ are of dimension $\eta$, then the expressions for $\omega$ and
$\tau$ are as follows. \bea \omega &=& 1+ \gamma {\bm \alpha} .({\bm  T}-\gamma {\bm  I}_{\eta})^{-1}. {\bm  e}_\eta , \lb{kagi_1} \\
\tau &=& 1- \lambda \gamma {\bm \alpha} .({\bm  T}-\gamma {\bm  I}_{\eta})^{-1}. {\bm  T}^{-1}{\bm 
e}_\eta . \lb{kagi_2_x} \eea {\bf  Proof:} Using the definition of $\omega$ and $\tau$, after little algebraic manipulation the the results (\ref{kagi_1}) and (\ref{kagi_2_x}) may be obtained. 
\lb{lt_1}\end{theorem}

\par In the following, we further define a few notations which are required to analyze
the queueing model under consideration. If we let  $A^{+}=A-\widehat{V}$ with
$A-\widehat{V}>0$, i.e., inter-arrival time is greater than remaining vacation time, then \bea F_{A^{+}}(x) & = &
\frac{\int_0^x\int_0^{\infty}f_{\widehat{V}}{(y)}f_{A}(y+s) \ dy \
ds}{Pr\{A>\widehat{V}\}} \nn \\ &=&
\frac{\int_0^x\int_0^{\infty}f_{\widehat{V}}{(y)}f_{A}(y+s) \ dy \
ds}{1-\omega},  \lb{kagi_6_x}\eea  where $A^{+}$ may be called excess inter-arrival time. Further, let us denote
$\widehat{A}^{+}$ and $\widetilde{A}^+$ as the remaining and elapsed
times of the excess inter-arrival time random variable $A^+$, respectively, then \bea
F_{\widehat{A}^{+}}(x)=F_{\widetilde{A}^{+}}(x)=\int_0^x\la_1(1-F_{A^{+}}(y))
\ dy,   \lb{gos_37}\eea where $\la_1=\int_0^{\infty} x \
dF_{A^{+}}(x)$ is the mean of the random variable $A^{+}.$ For an
important special case of phase-type inter-arrival time the
distribution of excess inter-arrival time $A^{+}$ is also phase-type which is proved in the
following theorem.
\begin{theorem}
If inter-arrival times follow a $PH$-distribution with
irreducible representation $({\bm \alpha},{\bm  T})$, where ${\bm \alpha}$ and ${\bm 
T}$ are of dimension $\eta$, then the distribution of $A^{+}$ is
also phase-type with representation $({\bm \alpha}_1,{\bm  T}_1)$, where
${\bm \alpha}_1$ and ${\bm  T}_1$ are of dimension $\eta$ and are given by
\bea {\bm \alpha}_1 &=&  \frac{\gamma}{1-\omega}{\bm \alpha} \Big(-{\bm  T}+\gamma {\bm 
I}_{\eta}\Big)^{-1}, \lb{kagi_14}
\\ {\bm  T}_1 &=& {\bm  T}. \lb{kagi_15}\eea
\nii {\bf  Proof:} Simple algebraic calculation will give the proof. 
\lb{lt_2}\end{theorem}

Further, we introduce a few more matrices which are required for the
analysis of this queueing model. With vacation ending, let ${\bm  V}_n~(n\geq 0)$ denote
the matrix of order $m\times m$ whose $(i,j)$th element represents
the conditional probability that $n$ customers are served during an
excess inter-arrival time $A^{+}$ and the service process passes to
phase $j$, provided at the initial instant of the previous arrival
epoch there were at least $(n+1)$ customers in the system and the server
was on a vacation with the service phase $i$. Then
 \bea {\bm  V}_n&=&\int_0^{\infty}{\bm 
P}(n,t)dF_{A^{+}}(t),~n\geq 0. \lb{kako_2}\eea Now, with vacation ending, let ${\bm 
V}_n^{*}$ denote the matrix of order $m\times m$ whose $(i,j)$th
element represents the probability that at least $n$ customers are served
during an excess inter-arrival period $A^{+}$ and the service
process is in phase $j$ with the server going on vacation at the end of the
$n$-th service completion, provided at the initial instant of
previous arrival epoch there were exactly $n$ customers in the
system and the server was on a vacation with the service phase $i$.
Then similar to the results derived above, we obtain \bea {\bm 
V}_n^{*} &=& \int_0^{\infty}\widetilde{{\bm  P}}(n,t)dF_{A^{+}}(t),
\quad n\geq 1. \lb{kako_3} \eea

\nii Further, with vacation ending, let ${\bm  \Delta}_n~(n\geq 0)$ denote the matrix of
order $m\times m$ whose $(i,j)$th element represents the limiting
probability that $n$ customers are served during an elapsed excess
inter-arrival time $\widetilde{A}^{+}$ with the service process
being in phase $j$, given that there were at least $(n+1)$ customers
in the system with the server being on a vacation with the service
phase $i$ at the beginning of the inter-arrival period. Then we have
the following expression for ${\bm  \Delta}_n$ \bea {\bm  \Delta}_n&=&
\int_0^{\infty}{\bm  P}(n,t)dF_{\widetilde{A}^{+}}(t)= \lambda_1
\int_0^{\infty}{\bm  P}(n,t)(1-F_{{A}^{+}}(t)) \ dt, \quad
n\geq 0. \lb{kako_4}\eea
\nii The relationships among the matrices ${\bm  V}_n$, ${\bm 
V}_n^{*}$ and ${\bm  \Delta}_n$ can be derived
as follows. \bea {\bm  \Delta}_0 &=& \lambda_1\Big({\bm  I}_m
-{\bm  V}_0\Big)(-{\bm  L}_0)^{-1}, \lb{kagi_7} \eea and \bea
{\bm  \Delta}_n  &=& \Big( {\bm  \Delta}_{n-1}{\bm  L}_{1}-\lambda_1
{\bm  V}_n\Big)(-{\bm  L}_0)^{-1}, \quad n\geq 1. \lb{kagi_8} \eea
Further, using \cite{cbp15}, it can be shown that  \bea {\bm 
V}_n^{*}&=& \frac{1}{\lambda_1}{\bm  \Delta}_{n-1}{\bm  L}_1,~n\geq 1.
\lb{newx_1}\eea
\nii Similarly, let ${\bm  \Delta}_n^{*}~(n\geq 1)$ denote the
$m\times m$ matrix whose $(i,j)$th element represents the limiting
probability that $n$ or more customers have been served during an
elapsed excess inter-arrival time $\widetilde{A}^{+}$ and the
service process is in phase $j$ with the server going on vacation at the
end of $n$th service completion, provided at the previous arrival
epoch the server was on a vacation with service phase $i$ and the
arrival lead the system to state $n$ or more customers. Then we
can write \bea {\bm  \Delta}_n^{*} &=& \int_0^{\infty}\widetilde{{\bm 
P}}(n,t)dF_{\widetilde{A}^{+}}(t)= \lambda_1
\int_0^{\infty}\widetilde{{\bm  P}}(n,t)(1-F_{\widetilde{A}^{+}}(t))
\ dt, \quad n\geq 1. \lb{kako_5} \eea
\nii
Now using the procedure discussed in \cite{cbp15}, one may derive the following relation:
 \bea {\bm  \Delta}_{n+1}^{*} &=&\left({\bm  \Delta}_n^{*}- {\bm  \Delta}_n
\right)\left(-{\bm  L}_0\right)^{-1}{\bm  L}_1,\quad n\geq
1.\label{phin}\eea

Finally, we define a few notations which are required to analyze
the queueing model under consideration. Let $A^{++}=A-(\widehat{V}+V)$
given that $A-(\widehat{V}+V)>0$. It is needless to mention that since $V$ is exponentially distributed, the distribution of $\widehat{V}+V$ will be Erlang of order two, which is a phase-type distribution with two states. Let the phase type representation of $\widehat{V}+V$ be denoted as $\bet=\left[\matrix{1.0 & 0.0}\right]$ with ${\bm 
U}=\left[\matrix{
  -\gamma & \gamma\cr
  0.0 & -\gamma}\right].$ To calculate the distribution function of $A^{++}$, we need to define $\omega_2$ which denotes the probability that $\widehat{V}+V$ exceeds
  an inter-arrival time $A$. Then,  \bea \omega_2 &=& \int_0^{\infty}
  Pr(A<\widehat{V}+V|\widehat{V}+V=x).f_{\widehat{V+V}}(x)
  \ dx \nn \\ &=& \int_0^{\infty} Pr(A<x).f_{\widehat{V}+V}(x) \ dx \\
  &=& \int_0^{\infty} F_A(x)f_{\widehat{V}+V}(x) \ dx. \lb{vis_3}\eea
   Similarly, if we let $\tau_2$ denote the probability that
  $\widehat{V}+V$ exceeds $\widehat{A}$, then \bea \tau_2 &=&
  \int_0^{\infty} F_{\widehat{A}}(x)f_{\widehat{V}+V}(x) \ dx.
  \lb{vis_4}\eea
Further, if an inter-arrival time is following $PH$ distribution with the above representation, then following the derivation as presented in Theorem \ref{lt_1}, we have
\bea  \omega_2 &=& 1+ (\frac{\gamma}{2})({\bm \alpha} \otimes \bet
).({\bm  T}\otimes {\bm  I}_{2}
+ {\bm  I}_{\eta}\otimes {\bm  U})^{-1}. ({\bm  e}_\eta \otimes {\bm  e}_{2}), \lb{kagi_1} \\
\tau_2 &=& 1- \lambda (\frac{\gamma}{2}) ({\bm \alpha} \otimes \bet ).({\bm  T}\otimes {\bm 
I}_{2}+ {\bm  I}_{\eta}\otimes {\bm  U})^{-1}. ({\bm  T}^{-1}{\bm 
e}_\eta \otimes {\bm  e}_{2}). \lb{kagi_2} \eea
Now the distribution function of $A^{++}$ may be derived as
\bea F_{A^{++}}(x) & = &
\frac{\int_0^x\int_0^{\infty}f_{\widehat{V}+V}{(y)}f_{A}(y+s) \ dy \
ds}{Pr\{A>\widehat{V}+V\}} \nn \\ &=&
\frac{\int_0^x\int_0^{\infty}f_{\widehat{V}+V}{(y)}f_{A}(y+s) \ dy \
ds}{1-\omega_2}.  \lb{kagi_6}\eea Similarly, let us denote
$\widehat{A}^{++}$ and $\widetilde{A}^{++}$ as the remaining and elapsed
times of the random variable $A^{++}$, then \bea
F_{\widehat{A}^{++}}(x)=F_{\widetilde{A}^{++}}(x)=\int_0^x\la_2(1-F_{A^{++}}(y))
\ dy,   \lb{gos_37}\eea where $\la_2=\int_0^{\infty} x \
dF_{A^{++}}(x)$ is the mean of the random variable $A^{++}.$ For an
important special case of phase-type inter-arrival time, the
distribution of $A^{++}$ is also phase-type whose representation can be obtained following Theorem 3.2 as $({\bm  {\bm \alpha}}_2,{\bm  T}_2)$, where
\bea {\bm \alpha}_2 &=&  \frac{\frac{\gamma}{2}}{1-\omega_2}({\bm \alpha} \otimes \bet
)\Big(-{\bm  T}\otimes {\bm  I}_{2}- {\bm  I}_{\eta}\otimes {\bm 
U}\Big)^{-1} \lb{kagi_14} \\ {\bm  T}_2 &=& {\bm  T} \otimes {\bm 
I}_{2}. \lb{kagi_15}\eea

Further, with another vacation ending, let ${\bm  C}_n~(n\geq 0)$ denote $m \times m$ matrix whose $(i,j)$th element represents
the conditional probability that $n$ customers are served during an
excess inter-arrival time $A^{++}$ and the service process passes to
phase $j$ with the server becoming idle (after completing service of the $n$-th customer), provided at the initial instant of the previous arrival
epoch there were at least $(n+1)$ customers in the system and the server
was on a vacation with the service phase $i$. Then
\bea {\bm  C}_n&=&\int_0^{\infty}{\bm 
P}(n,t)dF_{A^{++}}(t),~n\geq 0. \lb{kako_2_1}\eea

Similarly,  with the vacation ending, let ${\bm 
C}_n^{*}$ denote the matrix of order $m\times m$ whose $(i,j)$th
element represents the probability that at least $n$ customers are served
during an excess inter-arrival period $A^{++}$ and the service
process is in phase $j$ with the server becoming idle (after completing service of the $n$-th customer), provided at the initial instant of
previous arrival epoch there were exactly $n$ customers in the
system and the server was on a vacation with the service phase $i$.
Then, we can define \bea {\bm 
C}_n^{*} &=& \int_0^{\infty}\widetilde{{\bm  P}}(n,t)dF_{A^{++}}(t),
\quad n\geq 1. \lb{kako_3_1} \eea

Also, with the vacation ending, let ${\bm  \Phi}_n~(n\geq 0)$ denote $m \times m$ matrix whose $(i,j)$th element represents
the conditional probability that $n$ customers are served during an
elapsed excess inter-arrival time $\widetilde{A}^{++}$ and the service process passes to
phase $j$ with the server becoming idle (after completing service of the $n$-th customer), provided at the initial instant of the previous arrival
epoch there were at least $(n+1)$ customers in the system and the server
was on a vacation with the service phase $i$. Then, \bea {\bm  \Phi}_n&=&
\int_0^{\infty}{\bm  P}(n,t)dF_{\widetilde{A}^{++}}(t)= \lambda_2
\int_0^{\infty}{\bm  P}(n,t)(1-F_{{A}^{++}}(t)) \ dt, \quad
n\geq 0. \lb{kako_4_1}\eea

 Similarly, with the vacation ending, let ${\bm  \Phi}_n^{*}~(n\geq 1)$ denote the
$m\times m$ matrix whose $(i,j)$th element represents the limiting
probability that at least $n$ customers have been served during an
elapsed excess inter-arrival time $\widetilde{A}^{++}$ and the
service process is in phase $j$ with the server becoming idle (after completing service of the $n$-th customer), provided at the previous arrival
epoch the server was on a vacation with service phase $i$ with the
arrival leading the system to state $n$ customers. Then we
can write \bea {\bm  \Phi}_n^{*} &=& \int_0^{\infty}\widetilde{{\bm 
P}}(n,t)dF_{\widetilde{A}^{++}}(t)= \lambda_2
\int_0^{\infty}\widetilde{{\bm  P}}(n,t)(1-F_{{A}^{++}}(t))
\ dt, \quad n\geq 1. \lb{kako_5_1} \eea

\par  One may note here that the matrices ${\bm  V}_n~(n\geq 0)$ are required to obtain
other matrices by using the above relations
(\ref{kagi_7})-(\ref{kagi_8}).  
 The relationships among the matrices ${\bm  C}_n,~{\bm  C}_n^{*},~{\bm  \Phi}_n$ and ${\bm  \Phi}_n^{*}$ can be similarly obtained as for the  matrices ${\bm  V}_n,~{\bm  V}_n^{*},~{\bm  \Delta}_n$ and ${\bm  \Delta}_n^{*}$. These relationships are given below.
\bea {\bm  \Phi}_0 &=& \lambda_2\Big({\bm  I}_m
-{\bm  C}_0\Big)(-{\bm  L}_0)^{-1}, \lb{kagi_7_1} \\
{\bm  \Phi}_n  &=& \Big( {\bm  \Phi}_{n-1}{\bm  L}_{1}-\lambda_2
{\bm  C}_n\Big)(-{\bm  L}_0)^{-1}, \quad n\geq 1. \lb{kagi_8_1} \eea
\bea {\bm  \Phi}_1^{*}&=&
\left({\bm  I}_m - {\bm  \Phi}_0\right).\left(-{\bm 
L}_0\right)^{-1}{\bm  L}_1, \lb{xtan_2_1}\eea
\bea {\bm  \Phi}_{n+1}^{*} &=&\left({\bm  \Phi}_n^{*}- {\bm  \Phi}_n
\right)\left(-{\bm  L}_0\right)^{-1}{\bm  L}_1,\quad n\geq
1.\label{phin}\eea \bea {\bm 
C}_n^{*}&=& \frac{1}{\lambda_2}{\bm  \Phi}_{n-1}{\bm  L}_1,~n\geq 1,
\lb{newx_1_1}\eea The matrices ${\bm  C}_n$ are calculated
exactly the same way as we derive the matrices ${\bm  S}_n$, see Chaudhry et al. \cite{csa15}. For more information on $MSP$,
readers are referred to Bocharov \cite{nebb11}, Albores and Tajonar
\cite{nebb2} and Gupta and Banik \cite{bub122}.

\section{\bf  Analysis of $GI/MSP/1/\infty$ queue with single vacation}\lb{iid3}
We consider a $GI/MSP/1/\infty$ queueing system with single
vacation as described above. In the following subsections we obtain
steady-state distributions for this queueing system at different
epochs considering $\rho<1$.

\subsection{\bf  Stationary system-length distribution at pre-arrival epoch}\lb{iid4}
Consider the system just before arrival epochs which are taken as
embedded points. Let $t_0,t_1,t_2,...$~be the time epochs at which
arrivals occur and $t_k^{-}$ the time instant before $t_k$. The
inter-arrival times $T_{k+1}=t_{k+1}-t_k, ~k=0,1,2,\ldots$ are
i.i.d.r.v.'s with common distribution function $F_A(x)$. The state of
the system at $t_k^{-}$ is defined as
$\zeta_k=\{N_{t_k^{-}},J_{t_k^{-}}, \xi_{t_k^{-}}\}$ where
$N_{t_k^{-}}$ is the number of customers $n~(\geq 0)$ present in the
system including the one currently in service. Whereas
$J_{t_k^{-}}=\{j\},~1\leq j\leq m,$ denotes phase of the service
process and $\xi_{t_k^{-}}=l=0$ or 1 indicates that the server is on
vacation $(l=0)$ or busy $(l=1)$. In the limiting case, we define
the following probabilities: \bea \pi_{j,0}^{-}(n) &=&
\lim_{k\rightarrow
\infty}P\{N_{t_k^{-}}=n,~J_{t_k^{-}}=j,~\xi_{t_k^{-}}=0\},~ n\geq
0,~1\leq j\leq m, \nn \\ \pi_{j,1}^{-}(n) &=& \lim_{k\rightarrow
\infty}P\{N_{t_k^{-}}=n,~J_{t_k^{-}}=j,~\xi_{t_k^{-}}=1\},~ n\geq
0,~1\leq j\leq m, \nn \eea where $\pi_{j,0}^{-}(n)$ represents the
probability that there are $n~(\geq 0)$ customers in the system just
prior to an arrival epoch of a customer when the server is on
vacation with phase of the service process $j$. Similarly,
$\pi_{j,1}^{-}(n)$ denotes the probability that there are $n~(\geq
0)$ customers in the system just prior to an arrival epoch of a
customer when the server is in a busy (when $n\geq 1$) or dormant (when
$n=0$) state with phase of the service process $j$. Let
$\bpi_0^-(n)$ and $\bpi_{1}^-(n)$ be the row vectors of order
$1\times m$ whose $j$-th components are $\pi_{j,0}^{-}(n)$ and
$\pi_{j,1}^{-}(n),$ respectively.

\par Observing the state of the system at two consecutive embedded points, we
have an embedded Markov chain whose state space is equivalent to
$\Omega=\{(k,j,0),~k\geq 0, 1\leq j\leq m,\}\cup\{(n,j,1),~ n\geq
1,~1\leq j\leq m\}.$ Observing the system at two consecutive embedded Markov point, we have the following system of vector difference equations 
 \bea \bpi_0^-(0)&=& \sum_{k=0}^{\infty}\bpi_0^-(k)\Big((1-\omega){\bm  V}_{k+1}^{*}-{\bm  C}_{k+1}^{*}\Big)
+\sum_{n=0}^{\infty}\bpi_1^-(n)\Big( (1-\omega){\bm 
V}_{n+1}^{*}\Big),\label{c3}\\ \bpi_0^-(n)&=& \bpi_0^-(n-1)\omega {\bm  I}_m,\quad n\geq 1,\label{c4} \\
\bpi_1^-(0)&=& \sum_{k=0}^{\infty}\bpi_0^-(k){\bm 
C}_{k+1}^{*}+\sum_{n=0}^{\infty}\bpi_1^-(n)\Big(\omega ({\bm  V}_{n+1}^{*}+\sum_{j=0}^n{\bm  V}_j)+ (1-\omega)\sum_{j=0}^n{\bm  V}_j-\sum_{i=0}^n{\bm  S}_i\Big),\label{c3_idle}\\
\bpi_1^-(n)&=& \sum_{k=n-1}^{\infty}\bpi_0^-(k)(1-\omega){\bm 
V}_{k-n+1}+\sum_{j=n-1}^{\infty}\bpi_1^-(j){\bm  S}_{j-n+1} , \quad
n\geq 1.\lb{achhu_1} \eea \nii Multiplying (\ref{achhu_1})
 by $z^n$, summing from $n=1$ to $\infty$, after adding (\ref{c3_idle}) and using the vector-generating function
$\bpi_1^{-*}(z)=\sum_{n=0}^{\infty}\bpi_1^-(n)z^{n}$, we obtain \bea
\bpi_1^{-*}(z)[{\bm  I}_m- z {\bm  S}(z^{-1})]&=&
\sum_{j=0}^{\infty}\sum_{i=j+1}^{\infty} \bpi_0^-(j){\bm 
V}_{i}z^{j-i+1}\nn \\ &&+\sum_{k=0}^{\infty}\bpi_0^-(k){\bm 
C}_{k+1}^{*}+\sum_{n=0}^{\infty}\bpi_1^-(n)\Big(\omega ({\bm  V}_{n+1}^{*}+\sum_{j=0}^n{\bm  V}_j)\nn \\ && + (1-\omega)\sum_{j=0}^n{\bm  V}_j-\sum_{i=0}^n{\bm  S}_i\Big), \label{c5} \eea  leading to { \bea
\bpi_1^{-*}(z)&=&\frac{
\Big(\sum_{j=0}^{\infty}\sum_{i=j+1}^{\infty}
\bpi_0^-(j)(1-\omega){\bm 
V}_{i}z^{j-i+1}+{\bm  Y}\Big)Adj[{\bm  I}_m-z {\bm  S}(z^{-1})]}{det[{\bm  I}_m-z
{\bm  S}(z^{-1})]},   \label{c6} \eea \nii where ${\bm  Y}=\sum_{k=0}^{\infty}\bpi_0^-(k){\bm 
C}_{k+1}^{*}+\sum_{n=0}^{\infty}\bpi_1^-(n)\Big(\omega ({\bm  V}_{n+1}^{*}+\sum_{j=0}^n{\bm  V}_j)+ (1-\omega)\sum_{j=0}^n{\bm  V}_j-\sum_{i=0}^n{\bm  S}_i\Big).$ For further
analysis, we first determine an analytic expression for each
component of $\bpi_1^{-*}(z)$. Each component $\pi_{j,1}^{-*}(z)$
defined as $\pi_{j,1}^{-*}(z)=\sum_{n=0}^{\infty}\pi_{j,1}^-(n)z^n$
of the VGF $\bpi_1^{-*}(z)$ given in (\ref{c6}) being convergent in
$|z|\leq 1$ implies that $\bpi_1^{-*}(z)$ is convergent in $|z|\leq
1$. As each element of ${\bm  S}(z^{-1})$ is a rational function, see Chaudhry et al. \cite{cgg10}. Therefore, each element of $det[{\bm  I}_m-z {\bm  S}(z^{-1})]$
is also a rational function and we assume that \bee det[{\bm  I}_m-z
{\bm  S}(z^{-1})]=\frac{d(z)}{\varphi(z)}.\eee Equation (\ref{c6})
can be rewritten element-wise as \bea
\pi_{j,1}^{-*}(z)=\frac{\xi_j(z)}{d(z)},~~1\leq j\leq m,
\label{cs2}\eea where $\xi_j(z)$ is the $j$-th component of
$\Big(\sum_{j=0}^{\infty}\sum_{i=j+1}^{\infty}
\bpi_0^-(j)(1-\omega){\bm 
V}_{i}z^{j-i+1}+\sum_{k=0}^{\infty}\bpi_0^-(k){\bm 
C}_{k+1}^{*}\\+\sum_{n=0}^{\infty}\bpi_1^-(n)\Big(\omega ({\bm  V}_{n+1}^{*}+\sum_{j=0}^n{\bm  V}_j)+ (1-\omega)\sum_{j=0}^n{\bm  V}_j-\sum_{i=0}^n{\bm  S}_i\Big)\Big)Adj[{\bm  I}_m-z {\bm  S}(z^{-1})]\varphi(z)$. To
evaluate the vector in the numerator of equation (\ref{c6}), we show
that the equation $det[{\bm  I}_mz-{\bm  S}(z)]=0$ has exactly $m$
roots inside the unit circle $|z|=1$, see Chaudhry et al. \cite{csa12} Let these
roots be $\gamma_i~(1\leq i\leq m)$. Now, consider the zeros of the
function $d(z)$. Since the equation $det[{\bm  I}_mz-{\bm  S}(z)]=0$
has $m$ roots $\gamma_i$ inside the unit circle, the function
$det[{\bm  I}_m-z{\bm  S}(z^{-1})]$ has $m$ zeros $1/\gamma_i$ outside
the unit circle $|z|=1$. As $\pi_{j,1}^{-*}(z)$ is an analytic function of
$z$ for $|z|\leq 1$, applying the partial-fraction method, we
obtain\bea
\pi_{j,1}^{-*}(z)&=&\sum_{i=1}^{m}\frac{k_{ij}}{1-\gamma_iz},~1\leq
j\leq m,\label{c7}\eea where $k_{ij}$ are constants to be
determined. Now, collecting the coefficient of $z^n$ from both sides
of (\ref{c7}), we have\bea
\pi_{j,1}^-(n)&=&\sum_{i=1}^{m}k_{ij}\gamma^{n}_i,~1\leq j\leq
m,~n\geq 0. \label{c9}\eea Now we assume $\bpi_0^{-}(0)$ as \bea
\bpi_0^-(0)&=&\Bigg[b_1,b_2 ,\ldots,b_m\Bigg],\lb{kako_9}\eea where
$b_1,b_2 ,\ldots,b_m$ are $m$ arbitrary positive constants to be
computed as described below. Hereafter, we substitute
$\bpi_0^{-}(0)$ from Equation (\ref{kako_9}) into the Equation
(\ref{c4}) and obtain \bea \bpi_0^{-}(n)&=&
\bpi_0^{-}(0)\omega^n{\bm  I}_m,\quad n\geq 1. \lb{kagi_20} \eea From
(\ref{kako_9}) and (\ref{kagi_20}) we are able to express
$\bpi_0^{-}(n)~(n\geq 0)$ in terms of the $m$ constants $(b_1,b_2
,\ldots,b_m)$ and $\omega$ as defined above.

Next using (\ref{c9}) in (\ref{c3_idle}) and (\ref{achhu_1}) for
$n=0,2,\ldots,m-1$, we have \bea
\Bigg[\sum_{i=1}^{m}k_{i1},\sum_{i=1}^{m}k_{i2}
,\ldots,\sum_{i=1}^{m}k_{im}\Bigg] &=&
\sum_{k=0}^{\infty}\bpi_0^-(k){\bm 
C}_{k+1}^{*}
\nn \\ &&+\sum_{j=0}^{\infty}\Bigg[\sum_{i=1}^{m}k_{i1}\gamma_i^{j},\sum_{i=1}^{m}k_{i2}\gamma_i^{j}
,\ldots,\sum_{i=1}^{m}k_{im}\gamma_i^{j}\Bigg]\Big(\omega ({\bm  V}_{j+1}^{*}\nn \\ &&+\sum_{k=0}^j{\bm  V}_k)+ (1-\omega)\sum_{k=0}^j{\bm  V}_k-\sum_{i=0}^j{\bm  S}_i\Big),
\label{c10_nu_1} \\
\Bigg[\sum_{i=1}^{m}k_{i1}\gamma_i^{n},\sum_{i=1}^{m}k_{i2}\gamma_i^{n}
,\ldots,\sum_{i=1}^{m}k_{im}\gamma_i^{n}\Bigg] &=&
\sum_{k=n-1}^{\infty}\bpi_0^-(k)(1-\omega){\bm  V}_{k-n+1} \nn \\
&&+\sum_{j=n-1}^{\infty}\Bigg[\sum_{i=1}^{m}k_{i1}\gamma_i^{j},\sum_{i=1}^{m}k_{i2}\gamma_i^{j}
,\ldots,\sum_{i=1}^{m}k_{im}\gamma_i^{j}\Bigg]{\bm  S}_{j-n+1},  \lb{kako_10} \eea where $\bpi_0^{-}(n)~(n\geq 0)$ should be
used in terms of $\omega$ and the constants as given in Equations
(\ref{kako_9}) and (\ref{kagi_20}). Now Equations (\ref{c10_nu_1})
and (\ref{kako_10}) give $m^2$ simultaneous equations in $m(m+1)$
unknowns, $k_{ij}$'s $(1\leq i\leq m,~1\leq j\leq m)$ and
$b_i~(1\leq i\leq m)$. The other $m$ equations can be obtained
through equating the corresponding components of both sides of the
vector Equation (\ref{c3}) as follows. \bea \Bigg[b_1,b_2
,\ldots,b_m\Bigg]&=& \sum_{k=0}^{\infty}\bpi_0^-(k)\Big((1-\omega){\bm  V}_{k+1}^{*} -{\bm  C}_{k+1}^{*}\Big)\nn \\ &&
+\sum_{j=0}^{\infty}\Bigg[\sum_{i=1}^{m}k_{i1}\gamma_i^{j},\sum_{i=1}^{m}k_{i2}\gamma_i^{j}
,\ldots,\sum_{i=1}^{m}k_{im}\gamma_i^{j}\Bigg]\Big((1-\omega) {\bm 
V}_{j+1}^{*}\Big),
 \label{kagi_22} \eea
where $\bpi_0^{-}(k)~(k\geq 0)$ should be used in terms of $\omega$
and the constants as given in Equations (\ref{kako_9}) and
(\ref{kagi_20}). Finally, we have a total of $m^2+m=m(m+1)$
equations with $m(m+1)$ unknowns $k_{ij}$'s $(1\leq i\leq m,~1\leq
j\leq m)$ and $b_i~(1\leq i\leq m)$. One may note here that we
ignore any one component Equation of (\ref{kako_10}) for $n=m-1$ which
is a redundant equation and instead we use the normalization
condition given by \bea
\sum_{j=1}^m\pi_{j,1}^{-*}(1)+\sum_{n=0}^{\infty}\sum_{j=1}^m\pi_{j,0}^{-}(n)=1.
\lb{xz_1} \eea Above normalization condition can be simplified by
putting $z=1$ in Equation (\ref{c7}) leading to \bea
\sum_{j=1}^m\pi_{j,1}^{-*}(1)+\sum_{n=0}^{\infty}\sum_{j=1}^m\pi_{j,0}^{-}(n)
&=&\sum_{j=1}^m\sum_{i=1}^{m}\frac{k_{ij}}{1-\gamma_i}+\sum_{j=1}^m\frac{\pi_{j,0}^{-}(0)}{1-\omega}=1\label{c7_jen}\eea
Thus solving these $m(m+1)$ equations, we get $m(m+1)$ unknowns.


\subsection{\bf  Stationary system-length distribution at arbitrary epoch}\lb{iid6}
We now derive explicit expressions for the steady-state
queue-length distribution. Define
$\pi_{i,l}(n)=$ steady-state probability that $n~(\geq 0)$ customers are in the system at an arbitrary epoch with server busy $(l=1)$ or on vacation $(l=0)$ and the phase of the service process is $i~(1\le i\le m).$  In other words, $\bpi_l(n)=[\pi_{1,l}(n),\pi_{2,l}(n),\ldots,\pi_{m,l}(n)],~~n\geq
0;l=1$ or $n\geq 0;l=0,$  at an arbitrary epoch. Here $\bpi_l(n)$ is an arbitrary epoch stationary probability vector whose $j$-th component $\pi_{j,l}(n)~(1\leq j\leq m)$ is the steady-state probability that $n$ customers are in the system with server busy $(l=1)$ or vacation $(l=0).$ The classical
argument based on renewal theory relates the steady-state
system-length distribution at an arbitrary epoch to that at the
corresponding pre-arrival epoch. Using similar results of Markov
renewal theory and semi-Markov processes, see, e.g., \c{C}inlar
\cite{rev_19} or Lucantoni and Neuts \cite{rev_20}, we obtain\\ \bea
\bpi_0 (0)&=& \sum_{k=0}^{\infty}\bpi_0^-(k)\Big((1-\tau){\bm  \Delta}_{k+1}^{*}-{\bm  \Phi}_{k+1}^{*}\Big)
+\sum_{n=0}^{\infty}\bpi_1^-(n)\Big((1-\tau) {\bm 
\Delta}_{n+1}^{*}\Big), \lb{enst_1_nee5} \\
 \bpi_0 (n) & = &  \bpi_0^-(n-1)\tau {\bm  I}_m,
 \quad  n\geq 1, \lb{neast_2_nee6} \\ \bpi_1 (0)&=&  \sum_{k=0}^{\infty}\bpi_0^-(k){\bm 
 \Phi}_{k+1}^{*}+\sum_{n=0}^{\infty}\bpi_1^-(n)\Big(\tau (\Delta_{n+1}^{*}+\sum_{j=0}^n\Delta_j) + (1-\tau)\sum_{j=0}^n{\bm  \Delta}_j-\sum_{i=0}^n{\bm  \Omega}_i\Big),
\lb{enst_1_kh}
\\\bpi_1 (n) & = & \sum_{j=n-1}^{\infty}\bpi_0^-(j)(1-\tau){\bm 
\Delta}_{j-n+1}+ \sum_{j=n-1}^{\infty}\bpi^-(j){\bm  \Omega}_{j-n+1},
\quad n\geq 1. \lb{abdb13_nee7} \eea

\nii Note that since the service process is interrupted during
periods in which the server is on a vacation or in an idle state, it follows that
 \bea \sum\limits_{n=1}^{\infty}\bpi_1(n)\left({\bm 
L}_0+{\bm  L}_1\right)={\bm  0}, \lb{xtn_3}\eea which in turn implies
that \bea \sum\limits_{n=1}^{\infty}\bpi_1(n)=C\overline{\bpi},
\label{pibarc} \eea for some positive constant $C$. Thus, by post
multiplying the members of the previous equation by ${\bm  e}$, we
conclude that \bea \sum\limits_{n=1}^{\infty}\bpi_1(n){\bm  e}=C.
\label{pibarc_2} \eea Therefore, using (\ref{pibarc_2}) in
(\ref{pibarc}) we obtain \bea \frac{1}{\Big(
\sum_{n=1}^{\infty}\bpi_1(n){\bm 
e}\Big)}\sum\limits_{n=1}^{\infty}\bpi_1(n)=\overline{\bpi},
\lb{nod_3}\\ \mbox{i.e.,} ~\frac{1}{
\rho^{'}}\sum_{n=1}^{\infty}\bpi_1(n)=\overline{\bpi},
\lb{nod_33}\eea where $\rho^{'}=\sum_{n=1}^{\infty}\bpi_1(n){\bm  e}$
represents the probability that the server is busy. The above result
(\ref{nod_3}) is useful while performing numerical calculations.

\subsection{\bf  Queue-length distribution at post-departure epoch and their relation with pre-service epoch}\lb{iid6_2}
In this subsection, we derive the probabilities for the states of
the system immediately after a service completion takes place. Let
$\bpi^+(n)=[\pi^+_{1}(n),\pi^+_{2}(n),\ldots,\pi^+_{m}(n)],~n\geq
0,$ be the $1\times m$ vector whose $i$-th component $\pi^+_{i}(n)$
represents the post-departure epoch probability that there are $n$
customers in the queue immediately after a service completion of a
customer and the server is in phase $i$. The post-departure epoch
thus occurs immediately after the server has either reduced the
queue or become idle. Hence, using level-crossing arguments given in
Chaudhry and Templeton \cite[p. 299]{chau}, we have \bea
\bpi^+(n)&=&\frac{1}{\mu^{*}\rho^{'} }\bpi_1(n){\bm  L}_1,\quad n\geq
0.\label{postdeparture_2} \eea It may be noted that
$\sum_{n=1}^{\infty}\bpi_1(n){\bm  L}_1{\bm  e}=\mu^{*}\rho^{'}$, which represents the departure rate when the server is busy. \\

Let
$\bpi^{s-}(n)=[\pi^{s-}_{1}(n),\pi^{s-}_{2}(n),\ldots,\pi^{s-}_{m}(n)],~
n\geq 1,$ be the $1\times m$ vector whose $i$-th component
$\pi^{s-}_{i}(n)$ represents the pre-service epoch probability that
there are $n$ customers in the queue immediately before a service of
a customer takes place and the server is in phase $i$. The argument used
to find post-departure epoch probabilities may be based on the
distribution of the probabilities for the system at pre-service
epoch of a customer, the instant in time immediately before a real
service of a customer starts. Using the above arguments, we obtain
the following result. \bea  \bpi^{s-}(1)&=&
\bpi^+(1)+\bpi^+(0),\lb{visn_3}
\\ \bpi^{s-}(n)&=& \bpi^+(n),\quad n\geq 2. \lb{jht_1}\eea

\section{\bf  Performance measures}\lb{id99}
As state probabilities at various epochs are known, performance
measures can be easily obtained. The average number of customers in
the system (queue) at an arbitrary epoch are given by
$$L_s=\sum_{n=0}^{\infty}n \bpi_0 (n){\bm 
e}+\sum_{n=0}^{\infty}n \bpi_1 (n){\bm 
e},~~L_q=\sum_{n=1}^{\infty}(n-1) \bpi_0 (n){\bm 
e}+\sum_{n=1}^{\infty}(n-1) \bpi_1 (n){\bm  e}.$$
\subsection{\bf  Waiting-time analysis }\lb{id12}
In this section, we obtain the LST of waiting-time distribution of a
customer who is accepted in the system. Let $\bphi_k(\theta)$ be the
LST of the probability that $k$ customers will be served within a
time $x$ and the service process upon completion of service passes
to phase $j$, provided $k$ customers were in the system and the
service process was in phase $i$ at the beginning of service. Since
the probability that the service of a customer is completed in the
interval $(x,x+dx]$ is given by the matrix $e^{{\bm  L}_0x}{\bm 
L}_1dx$ and the total service time of $k$ customers is the sum of
their service times, ${\bphi}_1(s)$, the LST of service time with corresponding phase change, is given by  \bea
{\bphi}_1(s)=\int_{0}^\infty e^{-s x}e^{{\bm  L}_0x}{\bm 
L}_1dx=(s{\bm  I }_m-{\bm  L}_0)^{-1}{\bm 
L}_1,~\mbox{with}~\bphi_k(s)=\bphi_1^{k}(s), \quad k\geq 2.\lb{adb5}
\eea Further, we also need the LST of remaining vacation-time. It is
given by \bea f_{\widehat{V}}^*(s)&=& \frac{\gamma}{\gamma
+s}\lb{tomat_1} \eea
 Let $W_s^*(s)$ denote the LST of the actual waiting time
distribution of an arbitrary customer in system. The LST of the
waiting-time distribution in system can be derived as follows: \bea
~~ W_s^*(s)&=& \sum_{n=0}^{\infty}\bpi_1^-(n){\bphi}_1^{n+1}(s){\bm 
e}
+\sum_{n=0}^{\infty}\bpi_0^-(n)f_{\widehat{V}}^*(s).{\bphi}_1^{n+1}(s){\bm 
e}\nn \\ &=& \sum_{n=0}^{\infty}\bpi_1^-(n) \Big[(s{\bm  I }_m-{\bm 
L}_0)^{-1}{\bm  L}_1\Big]^{n+1}{\bm  e}_m \nn \\ &&
+\sum_{n=0}^{\infty}\bpi_0^-(n)\Big(\frac{\gamma}{\gamma
+s}\Big)\Big[(s{\bm  I }_m-{\bm  L}_0)^{-1}{\bm  L}_1\Big]^{n+1}{\bm 
e}_{m}, \nn \\ &=& \bpi_1^{-*}({\bphi}_1(s)){\bphi}_1(s)+(\frac{\gamma}{\gamma+s})\bpi_0^{-*}({\bphi}_1(s)){\bphi}_1(s), \lb{tomat_2}\eea where $bpi_0^{-*}(z)=\sum_{j=0}^{\infty}\bpi_0^{-}(j)z^j,~|z|\leq 1,$ defined similarly as $\bpi_1^{-*}$. One can find mean waiting time in the
system from Equation (\ref{tomat_2}) by differentiating it and
putting $s=0$. It is given by \bea W_s  &=&- W_s^{*(')}(0)\nn \\
&=& \sum_{n=0}^{\infty} \bpi_1^- (n)\sum_{j=0}^n(-{\bm  L}_0^{-1}{\bm 
L}_1)^{j}(-{\bm  L}_0^{-1})(-{\bm  L}_0^{-1}{\bm  L}_1)^{n-j}{\bm  e}_m+
\sum_{n=0}^{\infty}\bpi_0^-(n)\gamma (\gamma+s)^{-2}(-{\bm 
L}_0^{-1}{\bm  L}_1)^{n+1}{\bm  e}_m \nn \\ && +
\sum_{n=0}^{\infty}\bpi_0^-(n)\sum_{j=0}^n(-{\bm  L}_0^{-1}{\bm 
L}_1)^{j}(-{\bm  L}_0^{-1})(-{\bm  L}_0^{-1}{\bm  L}_1)^{n-j}{\bm  e}_m,
\nn \\ &=&  \sum_{n=0}^{\infty} \bpi_1^- (n)\sum_{j=0}^n(-{\bm 
L}_0^{-1}{\bm  L}_1)^{j}(-{\bm  L}_0^{-1}){\bm  e}_m+
\sum_{n=0}^{\infty}\bpi_0^-(n)(1/\gamma){\bm  e}_m \nn \\ && +
\sum_{n=0}^{\infty}\bpi_0^-(n)\sum_{j=0}^n(-{\bm  L}_0^{-1}{\bm 
L}_1)^{j}(-{\bm  L}_0^{-1}){\bm  e}_m. \quad [\mbox{As}~(-{\bm 
L}_0^{-1}{\bm  L}_1){\bm  e}_m={\bm  e}_m]\lb{adb7} \eea From the
Little's law, we can also get mean sojourn time as
$W_s(LL)=\frac{L_s}{\lambda }$ which may serve one of the
verifications while performing numerical computation.

\subsection{\bf  Expected length of busy and idle periods }\lb{id12_ap}

Since for this system, in the limiting case, the proportions of
times the server is busy and idle are $\rho^{'}$ and $1-\rho^{'}$,
respectively, we have
 \bea \frac{E(B)}{E(I)} &=& \frac{\rho^{'}}{1-\rho^{'}}, \lb{xton_1}\eea

\nii where $B$ and $I$ are random variables denoting the lengths of
busy and idle periods, respectively. We first discuss the mean busy
period $E(B)$, which is comparatively easy to evaluate. Let $N_q(t)$
denote the number of customers in system at time $t$ and $\xi_q(t)$ be the
state of the server, i.e., busy $(=1)$ or idle $(=0)$.
$\{N_q(t),\xi_q(t)\}$ enters the set of busy states, $\Upsilon
\equiv \{(0,1),(1,1),(2,1),\ldots\}$ at the termination of an idle
period. The conditional probability that $\{N_q(t),\xi_q(t)\}$
enters $(0,0),$ given that $\{N_q(t),\xi_q(t)\}$ enters $\Upsilon$,
is therefore $C\bpi^+(i){\bm  e},~i\geq 0,$ where
$C=\frac{1}{\bpi^+(0){\bm  e}}.$ Now $\{N_q(t),\xi_q(t)\}$ enters
$(i,1),~i\geq 0,$ irrespective of customers' arrival during a
service time, which may happen in expected time $E(S)$. Thus \bea
E(B)&=& \frac{\sum_{i=0}^{\infty}\bpi^+(i){\bm 
e}.E(S)}{\bpi^+(0){\bm  e}}=\frac{E(S)}{\bpi^+(0){\bm  e}}.
\lb{pan_11} \eea Substituting $E(B)$ from Equation (\ref{pan_11}) in
Equation (\ref{xton_1}), we obtain \bea E(I) &=&
\frac{1-\rho^{'}}{\rho^{'}}.\frac{E(S)}{\bpi^+(0){\bm 
e}}.\lb{xton_6}\eea
\section{\bf Approximation for system-length distributions based on one root}\lb{rd_2}

 In this context one may note that the heavy-traffic approximation was first investigated by J.F. C. Kingman (see \cite{kg61}) who showed that when the utilisation parameter ($\rho$) of an $M/M/1$ queue is near 1 a scaled version of the queue length process can be accurately approximated by a reflected Brownian motion. In this direction one can get several approximate results such as the tail
 probabilities of the queue-length distribution at a pre-arrival
 epoch, heavy- or light-traffic behaviour of the queue-length
 distributions based on the real root of $det[{\bm  I}_m-z{\bm 
 	S}(z^{-1})]=0$ which is closest to 1 and outside $|z|\leq 1$. The
existence of such a root has been discovered since long time in the
literature, e.g., Feller (\cite{fell_1}, pg. 276-277) calculated the
tail probabilities using a single root of the denominator which is
smallest root in absolute value. Also Chaudhry et al. \cite{chm90}
have given a formal proof of the existence of such a root. One may
note that it is not difficult to calculate this root numerically. In
this context it is worth mentioning that sometimes an approximate
value of this root may be used to get desired queue-length
distributions and this approximate value may be obtained in the
following way. We investigate the approximate root inside $|z|\leq
1$ by expanding the matrix of the left-hand side of the
characteristic equation $det[{\bm  I}_mz-{\bm  S}(z)]=0$ in powers of
$\rho$ as \bea {\bm  I}_mz-{\bm  S}(z)&=& {\bm  I}_mz-{\bm 
	S}(\rho+z-\rho) \nn \\ &=&{\bm  I}_mz- {\bm  S}(z-\rho)-\rho{\bm 
	S}^{'}(z-\rho)-\frac{\rho^2}{2!}{\bm  S}^{''}(z-\rho) +
o(\rho^2)\overline{{\bm  H}},\lb{rod_12} \eea where $\overline{{\bm 
		H}}$ is some unknown matrix and ${\bm  S}^{'}(.),~{\bm  S}^{''}(.)$
are the successive differentiation of ${\bm  S}(.)$ of order one and
two, respectively. Also one may note that in Equation
(\ref{rod_12}), as usual, $o(x)$ represents a function of $x$ with
the property that $\frac{o(x)}{x}\rightarrow 0$ as $x\rightarrow 0.$
Multiplying the right-hand side of the above Equation (\ref{rod_12})
by the vector $\overline{\bpi}$ from left and the vector ${\bm  e}$
from right, we may write the characteristic equation as follows \bea
z-\overline{\bpi}{\bm  S}(z-\rho){\bm  e}-\rho \overline{\bpi}{\bm 
	S}^{'}(z-\rho){\bm  e} - \frac{\rho^2}{2!}\overline{\bpi}{\bm 
	S}^{''}(z-\rho){\bm  e}+o(\rho^2)\overline{c}=0. \lb{kaki_1}\eea
where $\overline{c}=\overline{\bpi}\overline{{\bm  H}}{\bm  e}$ and is
a finite constant. Now an approximate value of the root for the
above described three cases is obtained as follows.
\begin{itemize}
	\item Light-traffic case: Applying $\rho\rightarrow 0+$ in Equation (\ref{kaki_1}) gives
	\bea z-\overline{\bpi}{\bm  S}(z-\rho){\bm  e}-\rho
	\overline{\bpi}{\bm  S}^{'}(z-\rho){\bm  e} -
	\frac{\rho^2}{2!}\overline{\bpi}{\bm  S}^{''}(z-\rho){\bm  e}=0,
	\lb{kaki_2}\eea which gives an approximate value of this root.
	\item Heavy-traffic case: Replacing $\rho$ by $(1-\rho)$ in Equation (\ref{kaki_1}) and applying the condition $\rho\rightarrow 1-$, we obtain
	\bea z-\overline{\bpi}{\bm  S}(z-1+\rho){\bm  e}-(1-\rho)
	\overline{\bpi}{\bm  S}^{'}(z-1+\rho){\bm  e} -
	\frac{(1-\rho)^2}{2!}\overline{\bpi}{\bm  S}^{''}(z-1+\rho){\bm  e}=0.
	\lb{kaki_3}\eea  Solving (\ref{kaki_3}) we get the desired value.
	\item Tail probabilities at a pre-arrival epoch: If $\rho_1$ denotes any arbitrary offered load numerically very close
	to $\rho$ then we replace $\rho$ by $(\rho-\rho_1)$ in Equation
	(\ref{kaki_1}) and applying the condition $\rho\rightarrow \rho_1$
	to get \bea z-\overline{\bpi}{\bm  S}(z-\rho+\rho_1){\bm 
		e}-(\rho-\rho_1) \overline{\bpi}{\bm  S}^{'}(z-\rho+\rho_1){\bm  e} -
	\frac{(\rho-\rho_1)^2}{2!}\overline{\bpi}{\bm 
		S}^{''}(z-\rho+\rho_1){\bm  e}=0. \lb{kaki_4}\eea \nii Hence, we can
	obtain the desired root, say $z_1$, by solving the Equation
	(\ref{kaki_4}) for $z$.
\end{itemize}

\nii Finally, it may be noted that once we obtain an approximate
value of a root, we can obtain the exact root through various
numerical methods. The equation $det[{\bm  I}_mz-{\bm  S}(z)]=0$ may
be used to find the original root which is closest to 1 in all the
above described cases. For the sake of completeness we present below
the procedure to calculate tail probabilities at a pre-arrival epoch
based on this one root.  To get the tail probabilities, assume \bea
\pi_{j,1}^{-}(n)\simeq k_{1,j}{z_1^n}=p_{j,1}^{a1}(n), \quad
n>n_{\epsilon},~1\leq j\leq m,\lb{rod_22}\eea where $n_{\epsilon}$
is chosen as the smallest integer such that
$|(\pi_{j,1}^{-}(n)-p_{j,1}^{a1}(n))/\pi_{j,1}^{-}(n)|<\epsilon,$
i.e., $|1-\frac{p_{j,1}^{a1}(n)}{\pi_{j,1}^{-}(n)}|<\epsilon .$ But,
since the probability $p_{j,1}^{a1}(n)$ follows a geometric
distribution with common ratio $z_1$, it is better to choose
$n_\epsilon$ such that
$|\frac{\pi_{j,1}^{-}(n)}{z_1\pi_{j,1}^{-}(n-1)}-1|<\epsilon.$ The
approximation gets better if more than one root, in ascending order
of magnitude, is used. It should however be mentioned that those
roots that occur in complex-conjugate pairs should be used in pairs.
Thus, the tail probabilities using three roots can be approximated
by \bea \pi_{j,1}^{-}(n)\simeq
\sum_{i=1}^3k_{i,j}{z_i^n}=p_{j,1}^{a3}(n), \quad
n>n_{\epsilon}^1,~1\leq j\leq m,\lb{rod_23}\eea where
$z_i~(i=1,2,3)$ are the roots in ascending order of magnitude and
$n_{\epsilon}^1$ may be chosen by
$|1-\frac{p_{j,1}^{a3}(n)}{\pi_{j,1}^{-}(n)}|<\epsilon,~n>n_{\epsilon}^1.$
Similar procedure may be adopted to calculate queue-length
distributions for the cases of light- and heavy-traffics. It may be remarked here that this root can
also be obtained accurately by simply using high precision of the software packages mentioned earlier. Similar way we can compute the waiting time distribution based on a few number of roots in case of light- and heavy-traffic.

\section{\bf   Numerical results and discussion}\lb{id15}

To demonstrate the applicability of the results obtained in the
previous sections, some numerical results have been presented in two
self explanatory tables. At the bottom of the tables, several
performance measures are given. Since various distributions can be
either represented or approximated by $PH$-distribution, we take
inter-arrival time distribution to be of $PH$-type having the
representation $({\bm \alpha}, {\bm  T})$, where ${\bm \alpha}$ and ${\bm 
T}$ are of dimension $\eta$. Then ${\bm  S}(z)$ can be derived as
follows using the procedure adopted in \cite{cgg10}.
\bea {\bm  S}(z)=({\bm  I}_m\otimes {\bm \alpha})({\bm  L}(z)\oplus{\bm  T})^{-1}({\bm 
I}_m\otimes{\bm  T}{\bm  e}_{\eta}), \lb{xxor_1} \eea with ${\bm  L}(z)\oplus{\bm 
T}=({\bm  L}(z)\otimes{\bm  I}_{\nu})+({\bm  I}_m\otimes{\bm  T})$, \
where $\oplus$ and $\otimes$ are used for Kronecker product and sum,
respectively. For the derivation of ${\bm  S}(z)$, see Chaudhry et al \cite{cgg10}. Knowing that each element of ${\bm  L}(z)$ is a
polynomial in $z$, each element of ${\bm  L}(z)\oplus{\bm  T}$ is also
a polynomial in $z$ and hence the determinant of $({\bm 
L}(z)\oplus{\bm  T})$ is a rational function in $z$. Thus, from
the above expression for ${\bm  S}(z)$, we can immediately say that
each element of ${\bm  S}(z)$ is a rational function in $z$ with the
same denominator. One may note that in case the degree of the polynomials in each element of ${\bm  S}(z)$ is very high, it may be difficult or time consuming to calculate the roots of the characteristic equation \bea det[{\bm  I}_mz-{\bm  S}(z)]=0. \lb{xxor_2} \eea This difficulty may be minimized by calculating ${\bm  S}(z)$ in the following way.
\bea {\bm 
S}(z)&=& \lim_{N \rightarrow \infty}\sum\limits_{n=0}^{N}{\bm 
S}_nz^n, \lb{xxor_3}\eea
where ${\bm 
S}_n$ may be obtained as proposed in Lucantoni \cite{abl2}.
 \\ We have carried out extensive numerical work
based on the procedure discussed in this paper by considering
different service matrices $MSP~({\bm  L}_0,{\bm  L}_1)$ and
phase-type inter-arrival time distribution $PH ({\bm \alpha}, {\bm 
T})$. All the calculations were performed on a PC having Intel(R)
Core 2 Duo processor @1.65 GHz with 8 GB  DDR2 RAM using MSPLE 18.
Further, though all the numerical results were carried out in high
precision, they are
reported here in 6 decimal places due to lack of space. \\
In Table 1, we have presented various epoch probabilities for a
$PH/MSP/1/\infty$ queue with exponential single vacation using our method described in this
paper. Vacation time is following exponential distribution with average number of vacations per unit of time is $\gamma=1.8$.
Inter-arrival time is $PH$-type and its representation is given by
$${\bm \alpha}=\left[\matrix{0.22& 0.33 & 0.45}\right]$$, $${\bm 
T}=\left[\matrix{
 -2.823 & 0.0 & 2.812\cr
  3.542 & -2.942 & 1.000 \cr
  1.710& 0.0 & -2.240}\right]$$  with $\la = 0.259558$.
  The $MSP$ matrices as\\ $${\bm  L}_0=\left[\matrix{
   -3.69939 & 0.01276 & 0.00572& 0.0 \cr
   0.01012 & -0.55759 & 0.0 &  0.00682\cr 0.0 & 0.02343 & -0.53152 & 0.48730 \cr 0.00649 & 0.55363 & 0.0 & -0.58531}\right],$$\\ $${\bm 
L}_1=\left[\matrix{
 3.65748 &  0.01727 &  0.0 & 0.00616\cr
 0.01353 & 0.00517& 0.52195& 0.0 \cr 0.00924& 0.0 & 0.0 & 0.01155 \cr 0.00561 &  0.0 & 0.00847 & 0.01111}\right]$$\\ with
  stationary mean service rate $\mu^*=1.121972$, lag-1 correlation coefficient $0.618173$ between successive service times and $\overline{\bpi}= \left[\matrix{
  0.264645 & 0.253046 & 0.254961 & 0.227348}\right]$ so that
  $\rho=\lambda/(\mu^*)=0.231341.$ To calculate system-length distribution we need
  to calculate the roots of \bea det[{\bm 
I}_mz-{\bm  S}(z)]=0,\lb{numa_2}\eea where $m=4$ as given above. Here
${\bm  S}(z)$ may be obtained by Equation (\ref{numa_1}) or (\ref{xxor_3}) for $N=70$, see Equation (\ref{xxor_1}). The $m =4$ roots of (\ref{numa_2}) inside $|z|<1$ are
evaluated. The corresponding $k_{ij}~(1\leq i\leq 4,~1\leq j\leq 4)$
and $b_i~(1\leq i\leq 4)$ values are calculated using the procedure
described in Section \ref{iid4}, see Appendix A. Now using Equation
(\ref{c9}), (\ref{kako_9}) and (\ref{kagi_20}), one can obtain
system-length distribution at pre-arrival epoch and after that using
relations (\ref{enst_1_nee5})-(\ref{abdb13_nee7}) the arbitrary
epoch probabilities may be derived, see Table 1.
\newpage
{\tiny\begin{center} {\bf Table 1:} System-length distributions at pre-arrival and arbitrary epoch.
\vspace{0.3cm}\\
 \begin{tabular}{|cccccc|} \hline
\mc{3}{|c}{Pre-arrival $\pi^{-}_{j,0}(n)$ } & \mc{3}{c|}{\& $\pi^{-}_{j,1}(n)$ }\\
$\pi^{-}_{j,0}(n)$ & $j=1$ & $j=2$ & $j=3$ & $j=4$ & $\sum_{j=1}^{m=4}$  \\
\hline 0 & 0.394155& 0.001935& 0.007826&  0.001078&  0.404994\\
1& 0.012922&  0.000063& 0.000256& 0.000035& 0.013277 \\
2& 0.000424&  0.000002&  0.000008& 0.000001&  0.000435 \\
3&  0.000014 &  0.000000&  0.000000 &  0.000000 & 0.000014\\
4&  0.000000 &  0.000000 &  0.000000 &  0.000000 &  0.000000 \\
5&   0.000000 &  0.000000 &  0.000000 &  0.000000 &  0.000000\\
$\vdots$& $\vdots$& $\vdots$ & $\vdots$ & $\vdots$ & $\vdots$   \\
\hline sum &  0.407514& 0.002001 &  0.008091& 0.001114 & 0.418720\\
\hline
$\pi^{-}_{j,1}(n)$ & $j=1$ & $j=2$ & $j=3$ & $j=4$ & $\sum_{j=1}^{m=4}$  \\
\hline 0& 0.418304 &  0.003215 &  0.006170 &  0.000057&  0.037466 \\
1&  0.013560&  0.007230&  0.009807&  0.006868&  0.037466  \\
2&  0.001028 &  0.007186&  0.007215&  0.006678 &  0.037466  \\
3&  0.000530 &  0.006045&  0.005157&  0.003092&  0.037466 \\
4&  0.000730 &  0.005106 &  0.004122 &  0.001808&  0.037466 \\
5& 0.000721& 0.004344&  0.003459 &  0.001355 &  0.037466 \\
$\vdots$& $\vdots$& $\vdots$ & $\vdots$ & $\vdots$ & $\vdots$   \\
\hline sum &  0.439280& 0.058473 &  0.056044&  0.027481& 0.581279 \\
\hline
\end{tabular}
\end{center}

\begin{center}
 \begin{tabular}{|cccccc|} \hline
\mc{3}{|c}{ Arbitrary $\pi_{j,0}(n)$ }& \mc{3}{c|}{\& $\pi_{j,1}(n)$} \\
 $\pi_{j,0}(n)$ & $j=1$ & $j=2$ & $j=3$ & $j=4$ & $\sum_{j=1}^{m=4}$   \\
\hline  0 &  0.304820&  0.001497 &  0.006065 &  0.000840 &  0.313221  \\
1&  0.057034&  0.000280&  0.001132&  0.000156 &  0.058602  \\
2 &  0.001870 &  0.000009&  0.000037&  0.000005&  0.001921 \\
3 &  0.000061&  0.000000&  0.000001&  0.000000&  0.000063 \\
4 &  0.000002&  0.000000&  0.000000& 0.000000 & 0.000002\\
5 & 0.000000 & 0.000000 & 0.000000 & 0.000000 & 0.000000  \\
$\vdots$& $\vdots$& $\vdots$ & $\vdots$ & $\vdots$ & $\vdots$   \\
\hline sum & 0.363788 &  0.001786&  0.007235&  0.001001 & 0.373810\\
\hline  $\pi_{j,1}(n)$ & $j=1$ & $j=2$ & $j=3$ & $j=4$ & $\sum_{j=1}^{m=4}$   \\
\hline  0&  0.445456 &  0.002265&  0.005750&  0.000927 &  0.454398 \\
1&  0.029024&  0.004327 &  0.007539 &  0.004676&  0.045566  \\
2 &  0.002520&  0.007486 &  0.007924&  0.007015 &  0.024945  \\
3 &   0.000626 &  0.006485 &  0.005700&  0.005660&  0.018471   \\
4 &   0.000504&  0.004673&  0.004313 &  0.003588 &  0.013079  \\
5 &  0.000457 &  0.003718&  0.003537&  0.002685 &  0.010397  \\
$\vdots$& $\vdots$& $\vdots$ & $\vdots$ & $\vdots$ & $\vdots$   \\
\hline sum &  0.481333&  0.050157 &  0.055174 &  0.039525& 0.626190  \\
\hline
$L_S$ & =1.020068, & $W_s$ & =3.955370, & $W_s(LL)$ & =3.930026.   \\
\hline
\end{tabular}
\end{center}}

\nii It may be noted that in the above numerical experiment, we can
find the conditional probability that the server is busy in phase
$i,~i=1,2.$ It is given by\\
$$\frac{1}{\rho^{'}}\sum_{n=1}^{\infty}\bpi_1(n)=\left[\matrix{
0.208843 & 0.278777 & 0.287695& 0.224685 }\right],$$ which matches with $\overline{\bpi}$
up to almost 2 decimal places. As shown in the above table, Little's law is satisfied up to two digits. These, to some extent, support the validity
of our analytical as well as numerical results.

\par  Next one may note here that this queueing model deals with
generally distributed inter-arrival time distribution. Theoretically,
it is possible to approximate any non-negative distribution
arbitrarily closely by a $PH$-type distribution, see Bobbio and
Telek \cite{Bobbio}, Bobbio et al. \cite{Bobbio et al.} and
references therein. Using these methods available in the literature,
one may estimate a general inter-arrival time distribution by a
$PH$-type distribution and calculate stationary distribution at
various epochs as stated above. As a demonstration, $PH$
approximation of inter-arrival time distributions which
are not phase-type, we consider the following example.
\par
Similar to
the above tables, in Table 2, we have presented various epoch probabilities for $LN/MAP/1/\infty$ queue with exponential
single vacation, where $LN$ stands for log-normal inter-arrival
time distribution. Vacation time is exponential  with the stationary mean vacation  rate
$\gamma=1.7$. Inter-arrival time is $LN$-type
and its probability density function and distribution functions are
given by $f_{A}(x)=\frac{1}{x\alpha\sqrt{2\pi}}e^{-\frac{(\ln
		(x)-\beta)^2}{2\alpha^2}},~\alpha=1.04,~\beta=0.215,~x> 0,$ and
$F_A(x)=\int_0^{x}f_{A}(u)du,~x>0,$ with $\lambda=0.469635$,
respectively.
The $MAP$ matrices have 4 phases in this case and their representation is given
by  $${\bm  L}_0=\left[\matrix{
	-2.69939 & 0.01276 & 0.00572& 0.0 \cr
	0.01012 & -0.55759 & 0.0 &  0.00682\cr 0.0 & 0.02343 & -1.53152 & 0.48730 \cr 0.00649 & 0.55363 & 0.0 & -0.58531}\right],$$ \\ $${\bm 
	L}_1=\left[\matrix{
	2.65748 &  0.01727 &  0.0 & 0.00616\cr
	0.01353 & 0.00517& 0.52195& 0.0 \cr 0.00924& 0.0 & 1.0 & 1.01155 \cr 0.00561 &  0.0 & 0.00847 & 0.01111}\right]$$ \\with
stationary mean service rate $\mu^*=1.112288$, lag-1 correlation coefficient 0.179796 between successive service times and $\overline{\bpi}= \left[\matrix{
	0.264645& 0.253046 &  0.254961& 0.227348}\right]$  so that
$\rho=\lambda/(\mu^*)=0.422224.$ . In the following we state a
procedure to discuss a $PH$ approximation of the $LN$-type
inter-arrival time. As discussed by Bobbio and
Telek \cite{Bobbio}, Bobbio et al. \cite{Bobbio et al.} and Pulungan
\cite{aph09}, we assume a minimal acyclic phase-type canonical
representation of the log-normal inter-arrival time distribution as follows:\\
$${\bm \alpha}=\left[\matrix{
	\alpha_1 & \alpha_2& \alpha_3}\right],~~$$ and $${\bm 
	T}=\left[\matrix{
	-\overline{t}_1& \overline{t}_1  & 0.0\cr
	0.0& -\overline{t}_2 & \overline{t}_2 \cr
	0.0 &0.0 & -\overline{t}_3}\right],$$ \\ where $\overline{t}_1,~\overline{t}_2,~
\overline{t}_3>0,$~$\overline{t}_1\leq \overline{t}_2\leq \overline{t}_3,$~$\alpha_1,~\alpha_2,~\alpha_3\geq
0,$ and ${\bm \alpha} {\bm  e}_{\eta}=1.0,~\eta=3.$ The probability density and
distribution function of this phase-type representation is given by
\bea f_{PH}(x)&=& {\bm \alpha} \ e^{{\bm  T}x} {{\bm  T}}^{0},\lb{plc}\\
F_{PH}(x)&=&1-  {\bm \alpha} \ e^{{\bm  T}x}{\bm  e}_{\eta}, \lb{plc_1}\eea
where ${{\bm  T}}^{0}$ is non-negative vector and satisfies
${\bm  T}{\bm  e}_{\eta}+{{\bm  T}}^{0}={\bm  0}$ and $\eta=3$. The
moments of this acyclic phase-type representation may obtained by
the following formulae: \bea \overline{\mu}_{i}^{PH}&=&
(-1)^{i}i!({\bm \alpha} {\bm  T}^{-i}{\bm  e}_{\eta}),\quad i\geq
1,~\eta=3.\eea Using the log-normal inter-arrival time probability density
function we can calculate moments by the following formulae: \bea
\overline{\mu}_{i}^{LN}&=& \int_0^{\infty}x^if_{A}(x)dx, \quad i\geq
1. \eea

The first four moments of the original log-normal distribution
have been targeted to match with corresponding moment of the acyclic
phase-type representation. This multi-objective problem can be
reduced to a single objective problem by assigning weights to
objectives. The non-linear programming (NLP) problem can be
formulated as:\\
\bea
\mbox{Minimize}     ~&~w_1 (\overline{\mu}_{1}^{LN}-\overline{\mu}_{1}^{PH})^2+w_2(\overline{\mu}_{2}^{LN}-\overline{\mu}_{2}^{PH})^2 + w_3(\overline{\mu}_{3}^{LN}-\overline{\mu}_{3}^{PH})^2 + w_4(\overline{\mu}_{4}^{LN}-\overline{\mu}_{4}^{PH})^2;\label{NLP_obj}\\
\mbox{Subject to}   ~&~ \alpha_i \geq 0,\quad i=1,~2,~3;\label{NLP_sub2}\\
~&~ {\bm \alpha} {\bm  e}_{\eta}=1.0,~\eta=3;\label{NLP_sub3}\\                    ~&~ \overline{t}_i > 0,\quad i=1,~2,~3;
\label{NLP_2}\eea  \nii where $w_i$'s ($i=1,~2,~3,~4$) are the weights.
We always set $w_i=1,~i=1,~2,~3,~4$ except in cases when numerical
values of moments are very high. For example when numerical value of
fourth moment is very high we set $w_4=0$ and NLP problem gives
better result by assigning $w_4=0$ in lieu of $w_4=1.$ One may note
that after solving the above NLP we can completely specify an
approximate acyclic phase-type representation of the given weibull
distribution. For better approximation one may increase the order of
acyclic phase-type distribution and in that case obtaining a
solution to the above NLP may pose some problem due to increase in the
number of variables. To overcome this problem we supply eigenvalues
of the matrix ${\bm  T}$, i.e., $\overline{t}_1$ $\overline{t}_2$ and $\overline{t}_3 $ using the following
procedure. We calculate 20 moments from $f_{A}(x)$ to construct an
approximate $f_{A}^*(s)$ through the Pad\'{e} approximation
[2/3] as: \bea f_{A}^*(\theta) \simeq
\frac{1.0+64.325780
	s+479.132648s^2}{1.0+66.455094s+613.950181s^2+904.240262s^3}.
\lb{tgh_1_2}\eea where
$[2/3]$ stands for a rational function with degree of numerator
polynomial 2 and degree of denominator polynomial 3 in the variable $s$. Now after
equating denominator of (\ref{tgh_1_2}) equal to zero we obtain three
roots as $-0.017943,~-0.112322$ and $-0.548702.$ Hereafter, we assign
eigenvalues of the matrix ${\bm  T}$ as $\overline{t}_1=0.017943,~\overline{t}_2=0.112322$ and
$\overline{t}_3=0.548702.$ After this, our job is to solve the above NLP and obtain
the vector ${\bm \alpha}.$ Solution of above NLP with assignment of weights
as $w_i=1,~i=1,~2,~3,~4$ gives $\alpha_1=0.000024,~\alpha_2=0.034291$ and
$\alpha_3=0.965685.$ Therefore, we are able to specify the vector
${\bm \alpha}$ as well as the matrix ${\bm  T}$ corresponding to an
approximate phase-type distribution of the log-normal inter-arrival time
distribution. We calculate $\lambda$ using these ${\bm \alpha}$ and
${\bm  T}$ values and found $\lambda=0.469635$ which is exactly the
same as the one calculated using the original log-normal density function given above. In
the following we present graphs of the density and distribution
functions (see Figure 1 and 2, respectively) for the original
log-normal inter-arrival time and the corresponding phase-type
approximation.
\begin{figure}[h]
	\begin{center}
		\includegraphics[height=6cm,scale=0.5]{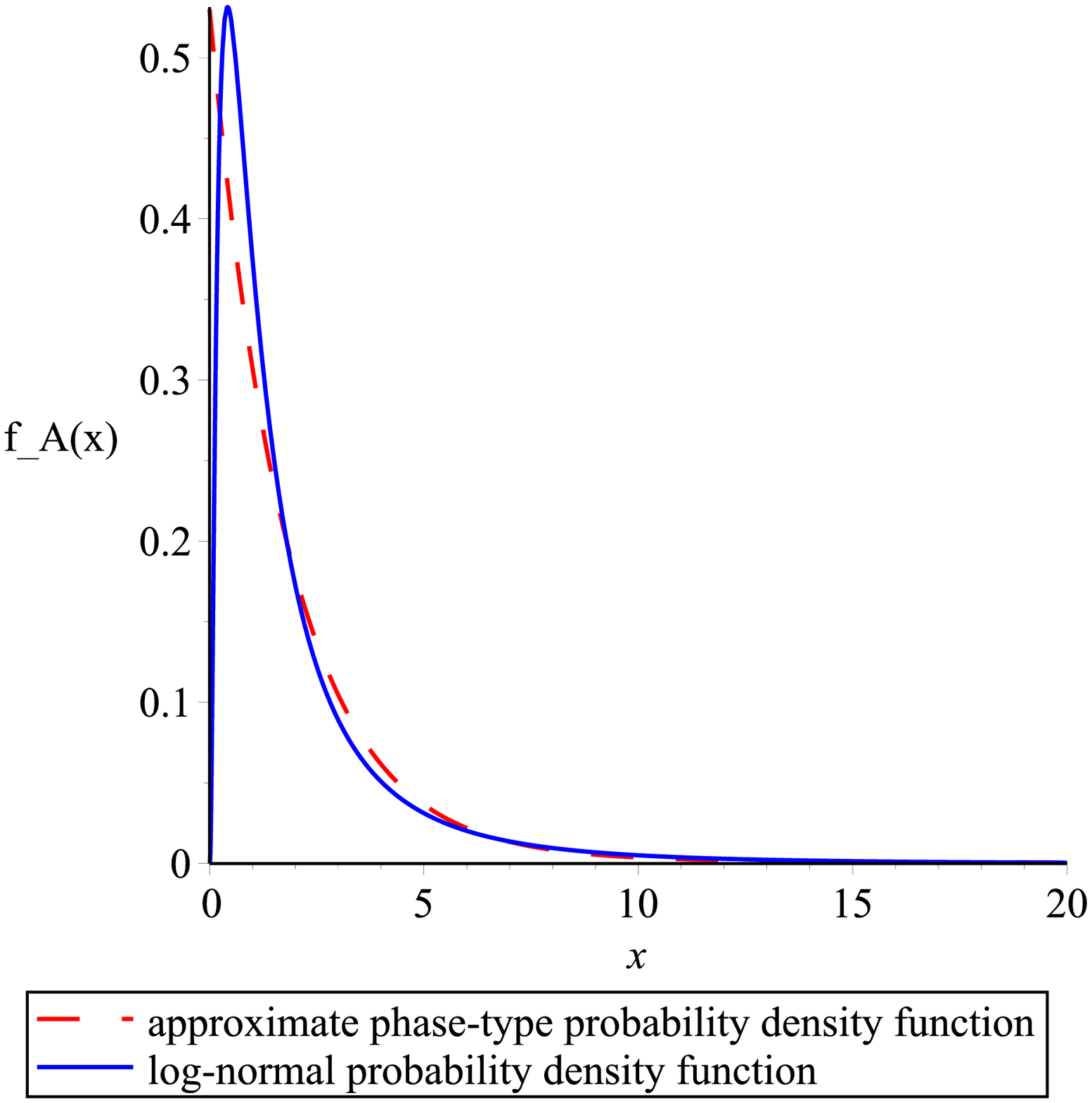} \caption{
			$x$ versus $f_A(x)$} \end{center}\end{figure}\lb{fg_1}
\begin{figure}[h]
	\begin{center}
		\includegraphics[height=6cm,scale=0.5]{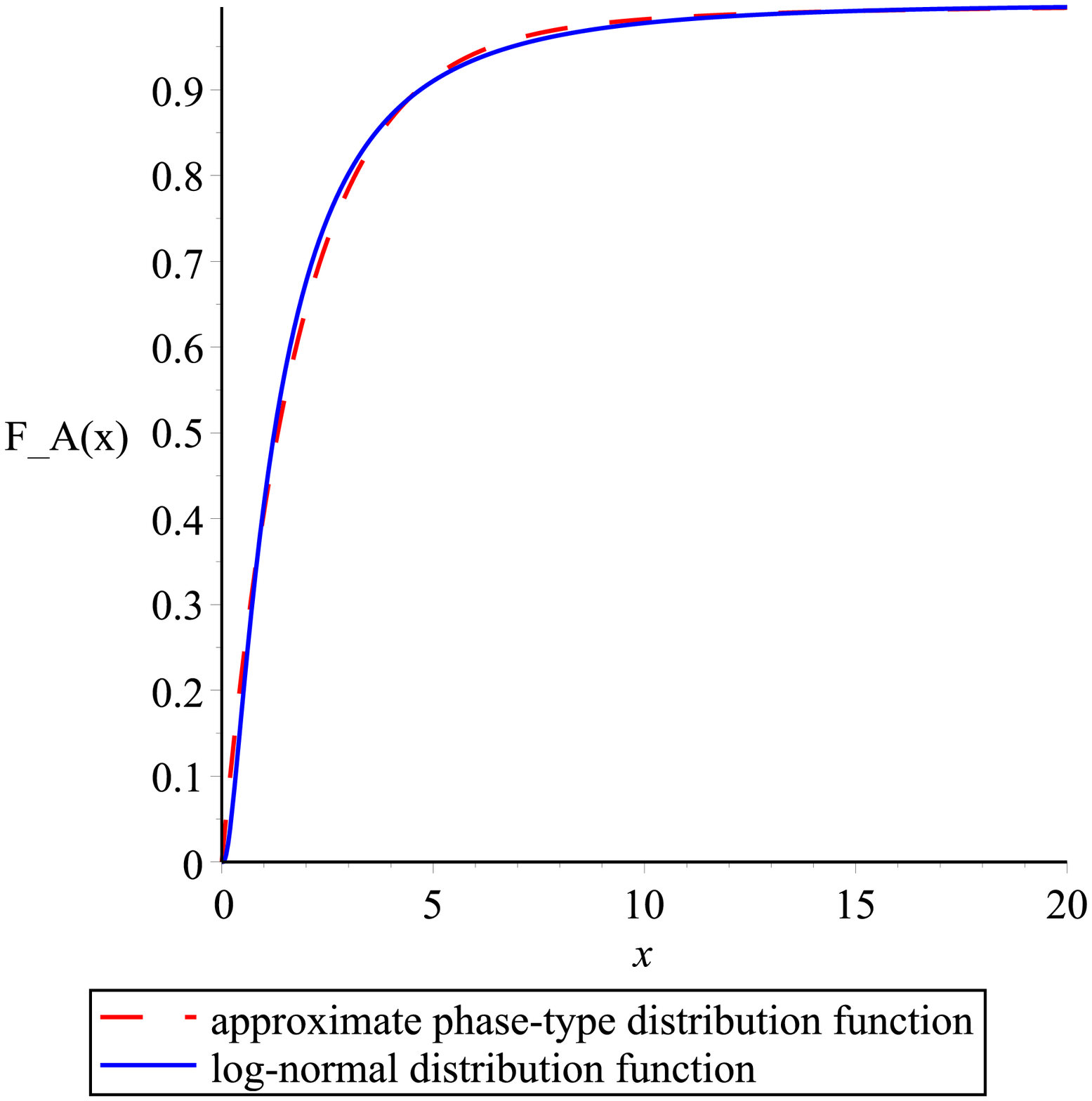} \caption{
			$x$ versus $F_{A}(x)$} \end{center}\end{figure}\lb{fg_2}
\par  Now we consider these approximate phase-type
representations for the corresponding inter-arrival and vacation
time distributions and using similar procedure as described
in Table 1, we have computed stationary system length
distribution at pre-arrival and arbitrary epochs, see Table 2.
\newpage
{\tiny \begin{center} {\bf Table 2:} System-length distributions at pre-arrival and arbitrary epoch.\vspace{0.3cm}\\
 \begin{tabular}{|cccccc|} \hline
\mc{3}{|c}{Pre-arrival $\pi^{-}_{j,0}(n)$ } & \mc{3}{c|}{\& $\pi^{-}_{j,1}(n)$ }\\
$\pi^{-}_{j,0}(n)$ & $j=1$ & $j=2$ & $j=3$ & $j=4$ & $\sum_{j=1}^{m=4}$  \\
\hline 0 &  0.190353&  0.001284 &  0.028939 &  0.000850 &  0.221426\\
1&  0.044952&  0.000303&  0.006834&  0.000201 &  0.052291\\
2&  0.010616 &  0.000072 &  0.001614 &  0.000047 &  0.012349  \\
3&  0.002507 &  0.000017 &  0.000381&  0.000011&  0.002916  \\
4&   0.000592&  0.000004 &  0.000090&  0.000003 &  0.000689 \\
5& 0.000140 &  0.000001&  0.000021 & 0.000001&  0.000163\\
$\vdots$& $\vdots$& $\vdots$ & $\vdots$ & $\vdots$ & $\vdots$   \\
\hline sum &   0.249203&  0.001680&  0.037886 &  0.001113 & 0.289883\\
\hline
$\pi^{-}_{j,1}(n)$ & $j=1$ & $j=2$ & $j=3$ & $j=4$ & $\sum_{j=1}^{m=4}$  \\
\hline 0& 0.268269 &  0.149411 &  0.030733 &  0.000812&  0.136142 \\ 1&  0.090191&  0.009359&  0.024661&  0.011930&  0.136142 \\
2&  0.027599&  0.012151 &  0.016200&  0.013089&  0.136142\\
3&  0.008699&  0.012077&  0.011171 &  0.009973&  0.136142 \\
4&  0.003498 &  0.010862 & 0.008328&  0.007005&  0.136142\\
5& 0.002095 &  0.009355 &  0.006560&  0.005010 &  0.136142\\
$\vdots$& $\vdots$& $\vdots$ & $\vdots$ & $\vdots$ & $\vdots$   \\
\hline sum & 0.410444 &  0.102596&  0.128800&  0.068276& 0.710117 \\
\hline
\end{tabular}
\end{center}
\begin{center}
 \begin{tabular}{|cccccc|} \hline
\mc{3}{|c}{ Arbitrary $\pi_{j,0}(n)$ }& \mc{3}{c|}{\& $\pi_{j,1}(n)$} \\
 $\pi_{j,0}(n)$ & $j=1$ & $j=2$ & $j=3$ & $j=4$ & $\sum_{j=1}^{m=4}$   \\
\hline  0 &  0.224084&  0.001540 &  0.042695&   0.001126 &  0.269445  \\
1&  0.042409&  0.000286 &  0.006447&  0.000189&  0.049332\\
2 & 0.010015&  0.000067&  0.001523 &  0.000045 &  0.011650\\
3 &  0.002365&  0.000016 &  0.000360 &  0.000011 &  0.002751  \\
4 &  0.000558&  0.000004 &  0.000085 &  0.000002&  0.000650\\
5 & 0.000132 &  0.000001&  0.000020&  0.000000 &  0.000153 \\
$\vdots$& $\vdots$& $\vdots$ & $\vdots$ & $\vdots$ & $\vdots$   \\
\hline sum &  0.279605&  0.001915 &  0.051135&  0.001374  & 0.334030 \\
\hline  $\pi_{j,1}(n)$ & $j=1$ & $j=2$ & $j=3$ & $j=4$ & $\sum_{j=1}^{m=4}$   \\
\hline  0 &  0.283429&  0.000799 &  0.035623&  0.000270 &  0.320121   \\ 1 &  0.044144&  0.006548&  0.016608 &  0.008447&  0.075748  \\
2 & 0.026168 &  0.012352 &  0.015941&  0.013083 &  0.067544 \\
3 &  0.008292 &  0.012130&  0.011055 &  0.011245&  0.042723 \\
4 &
 0.003210 &  0.010447&  0.008253&  0.008496 &  0.030406 \\
5 &  0.001761&  0.008642 &  0.006494&  0.006310&  0.023208\\
$\vdots$& $\vdots$& $\vdots$ & $\vdots$ & $\vdots$ & $\vdots$   \\
\hline sum & 0.374541 &  0.092439 & 0.124715&  0.074275& 0.665970 \\
\hline
$L_S$ & =1.817167, & $W_s$ & = 3.776260,& $W_s(LL)$ & = 3.869320.  \\
\hline
\end{tabular}
\end{center}}

\nii It may be noted that in the above numerical experiment, we can
find the conditional probability that the server is busy in phase
$i,~i=1,2.$ It is given by
$$\frac{1}{\rho^{'}}\sum_{n=1}^{\infty}\bpi_1(n)=\left[\matrix{0.263443& 0.264971 & 0.257603 &  0.213982 }\right],$$ which is $\overline{\bpi}$ up to
almost 2 decimal places as was anticipated. Also, Little's law is satisfied up to two digits after rounding off. These facts confirm our analytical as well as numerical results.



\section{Conclusions and future scope}\lb{anew}

In this paper, we have successfully analyzed the
$PH/MSP/1/\infty$ queue with single exponential vacation. We have suggested a procedure to obtain the steady-state
distributions of the number of customers in the system at
pre-arrival, arbitrary and post-departure epochs. Similar kind of analysis that is described in this paper may work for the corresponding queueing systems under multiple exponential vacations of the server or for batch arrival or batch service queues,
i.e., $PH^{[X]}/MSP/1/\infty$ queue or $PH/MSP^{(a,b)}/1/\infty$ queue
with exponential single or multiple vacations. One may be interested in analyzing the same queueing model with different type of
vacation policies, e.g., multiple adaptive
vacation(s) and working vacation(s). Another area of interest may be to find the approximations for the tail of the waiting-time distribution a swell as an approximation for the waiting-time distribution in cases of heavy- and light-traffic. These problems are left for future investigations.

\nii {\bf    Appendix A}

\nii The roots used in Table 1 are given as follows:
$\gamma_1=-0.010443,~\gamma_2=0.853818,~\gamma_3=0.160075-0.064040\mathrm{i},~\gamma_4=0.160075+0.064040\mathrm{i}.$ The corresponding
$k_{ij}~(1\leq i\leq 4,~1\leq j\leq 4)$ and $b_i~(1\leq i\leq 4)$
values are as follows: $k_{1,1}=0.250531,~k_{1,2} =0.009814,~k_{1,3} =0.017357,~k_{1,4} = 0.093024,~
k_{2,1}=0.001666, ~k_{2,2} =0.009563 ,~k_{2,3}=0.007588 ,~k_{2,4}=0.002873,~k_{3,1}=0.083053+0.092411\mathrm{i}, ~k_{3,2}= -0.008081-0.013702\mathrm{i},~k_{3,3}=-0.009387-0.050870\mathrm{i} ,~ k_{3,4}=-0.047920-0.161839\mathrm{i} ,~k_{4,1}=0.083053-0.092411\mathrm{i} ,~k_{4,2}=-0.008081+0.013702\mathrm{i} ,~k_{4,3}=-0.009387+0.050870\mathrm{i} ,~ k_{4,4}=-0.047920+0.161839\mathrm{i}$ and
$b_1 =0.394155,~ b_2 = 0.001935,~ b_3 = 0.007826,~ b4 =0.001078$.\\
The roots used in Table 2 are given below.
$\gamma_1=0.194296,~\gamma_2=0.831411,~\gamma_3=0.328896-0.025835\mathrm{i}$ and $\gamma_4=0.328896+0.025835\mathrm{i}.$ The
corresponding $k_{ij}~(1\leq i\leq 4,~1\leq j\leq 4)$ and
$b_i~(1\leq i\leq 4)$ values are as follows:
$k_{1,1}=-0.205961,~k_{1,2} =0.001301 ,~k_{1,3} =-0.036359 ,~k_{1,4} = 0.147868,~
k_{2,1} =0.005287,~k_{2,2} =0.024221,~k_{2,3}=0.015824,~k_{2,4}=0.010081 ,~k_{3,1}= 0.234467 -0.549961\mathrm{i}, ~k_{3,2}=-0.012014 -0.060543\mathrm{i} ,~k_{3,3}=0.025634+0.033050\mathrm{i} ,~ k_{3,4}=-0.078568+0.512859\mathrm{i},~k_{4,1}= 0.234467 +0.549961\mathrm{i} ,~k_{4,2}=-0.012014 +0.060543\mathrm{i} ,~k_{4,3}=0.025634-0.033050\mathrm{i},~ k_{4,4}=-0.078568-0.512859\mathrm{i}$ and
$b_1 = 0.190353,~ b_2 = 0.001284,~ b_3 = 0.028939,~ b_4 =0.000851$.

{\bf  Acknowledgement}\lb{id155}
 The second author
 was supported partially by NSERC under research grant number RGPIN-2014-06604.

\end{document}